\documentclass[reqno,10pt,centertags,draft]{amsart}
\usepackage{amsmath,amsthm,amscd,amssymb,latexsym,upref}

\newcommand{\bbN}{{\mathbb{N}}}
\newcommand{\bbR}{{\mathbb{R}}}

\newcommand{\bbZ}{{\mathbb{Z}}}
\newcommand{\bbC}{{\mathbb{C}}}

\newcommand{\calA}{{\mathcal A}}
\newcommand{\calB}{{\mathcal B}}

\newcommand{\calD}{{\mathcal D}}

\newcommand{\calS}{{\mathcal S}}

\newcommand{\dott}{\,\cdot\,}
\newcommand{\no}{\notag}
\newcommand{\lb}{\label}
\newcommand{\f}{\frac}

\newcommand{\wti}{\widetilde}
\newcommand{\ta}{\wti{\alpha}}
\newcommand{\tb}{\wti{\beta}}
\newcommand{\tg}{\wti{\gamma}}

\newcommand{\rank}{\text{\rm{rank}}}

\newcommand{\bi}{\bibitem}
\newcommand{\hatt}{\widehat}
\newcommand{\beq}{\begin{equation}}
\newcommand{\eeq}{\end{equation}}
\newcommand{\ba}{\begin{align}}
\newcommand{\ea}{\end{align}}
\newcommand{\vp}{\varphi}
\newcommand{\vt}{\vartheta}

\renewcommand{\Re}{\text{\rm
Re}}
\renewcommand{\Im}{\text{\rm Im}}
\renewcommand{\ln}{\text{\rm
ln}}
\renewcommand{\ge}{\geqslant}
\renewcommand{\le}{\leqslant}

\allowdisplaybreaks
\numberwithin{equation}{section}

\newtheorem{theorem}{Theorem}[section]
\newtheorem{lemma}[theorem]{Lemma}
\newtheorem{corollary}[theorem]{Corollary}
\newtheorem{remark}[theorem]{Remark}

\theoremstyle{definition}
\newtheorem{definition}[theorem]{Definition}
\newtheorem{hypothesis}[theorem]{Hypothesis}

\begin{document}

\title[Weyl--Titchmarsh theory for singular finite difference systems]{
On Weyl--Titchmarsh Theory for Singular Finite Difference Hamiltonian
Systems}
\author[S.\ Clark and F.\ Gesztesy]{Steve Clark and Fritz Gesztesy}
\address{Department of Mathematics and Statistics, University
of Missouri-Rolla, Rolla, MO 65409, USA}
\email{sclark@umr.edu}
\urladdr{http://www.umr.edu/\~{ }clark}
\address{Department of Mathematics,
University of
Missouri, Columbia, MO
65211, USA}
\email{fritz@math.missouri.edu\newline
\indent{\it URL:}
http://www.math.missouri.edu/people/fgesztesy.html}
\dedicatory{Dedicated with great pleasure to Norrie Everitt on the
occasion of his 80th birthday.}
\subjclass{Primary 34B20, 39A10;  Secondary 34B27, 39A70}
\keywords{Discrete canonical systems, difference operators,
Weyl--Titchmarsh theory, Green's functions}
\thanks{To appear in {\it J. Comput. Appl. Math.}}

\begin{abstract}
We develop the basic theory of matrix-valued
Weyl--Titchmarsh M-functions and the associated Green's matrices for
whole-line and half-line self-adjoint Hamiltonian finite
difference systems with separated boundary conditions. 
\end{abstract}

\maketitle

\section{Introduction}\label{s1}

This paper can be viewed as a natural continuation of our recent work
on matrix-valued Schr\"odinger and Dirac-type operators (cf.\ \cite{CG01}
and \cite{CG02}) to discrete Hamiltonian systems (i.e., Hamiltonian
systems of difference equations). These investigations are part of a
larger program which includes the following: \\
$(i)$ A systematic asymptotic expansion of Weyl--Titchmarsh matrices and
Green's matrices as the spectral parameter tends to infinity
(\cite{CG01}, \cite{CG02}). \\
$(ii)$ The derivation of trace formulas for such systems (\cite{CG02},
\cite{CGHL00}, \cite{GH97}). \\
$(iii)$ The proof of certain uniqueness theorems (including Borg and
Hochstadt-type theorems) for the operators in question (\cite{BGMS03},
\cite{CG02}, \cite{CGHL00}, \cite{GKM02}, \cite{GS00}). \\
$(iv)$ The application of these results to related integrable systems
(cf.\ \cite{BGMS03}, \cite{GS03}).

Before we describe the content of this paper in more detail, it is
appropriate to briefly comment on the literature devoted to general
$2m\times 2m$ Hamiltonian systems ($m\geq 2$) and their (inverse)
spectral theory as it relates to the topics of this paper and the next
one in our series (see \cite{CGR04}). Due to the enormous amount of
interest generated by continuous Hamiltonian systems over the past twenty
years, we are forced to focus primarily on references in connection with
discrete Hamiltonian systems, but we refer the reader to
\cite{BGMS03}, \cite{CG01}--\cite{CGHL00}, \cite{GKM02}, and \cite{GS03}
which provide  extensive documentation of pertinent material. The basic
Weyl--Titchmarsh theory of regular Hamiltonian systems can be found in
Atkinson's monograph \cite{At64}; Weyl--Titchmarsh theory of singular
Hamiltonian systems and their basic spectral theory was developed by
Hinton and Shaw and many others (see, e.g., \cite[Sect.\ 10.7]{Hi69},
\cite{HS93}--\cite{KR74}, \cite{Kr89a}--\cite{LM02}, \cite{Or76},
\cite{Ro60}--\cite{Sa94a}, \cite[Ch.~9]{Sa99a}, \cite[]{We87}
and the references therein); the corresponding theory for Jacobi
operators can be found in \cite[Sect.\ VII.2]{Be68}, \cite{Fu76},
\cite[Ch.\ 10]{Sa97} and the literature therein. Deficiency indices of
matrix-valued Jacobi operators are studied in \cite{KM98}--\cite{KM01}.
Inverse spectral and scattering theory for matrix-valued finite
difference systems and its intimate connection to matrix-valued
orthogonal polynomials and the moment problem are treated in \cite{AG94},
\cite{AN84}, \cite[Sect.\ VII.2]{Be68}, \cite{DL96}--\cite{DV95},
\cite{Ge82}, \cite{Lo99}, \cite{MBO92}, \cite{Os00}, \cite{Os02},
\cite[Ch.\ 8]{Sa97}, \cite{Se80}. Finally, connections with  nonabelian
completely integrable systems are discussed in \cite{BGS86},
\cite{BG90}, \cite{Os97}, \cite[Chs.\ 9, 10]{Sa97}.

In spite of these activities, the reader might perhaps be surprised to
hear that Weyl--Titchmarsh theory for general discrete Hamiltonian systems
appears to be underdeveloped. The only notable exceptions to this
statement of course being the special case of matrix-valued Jacobi
operators which are described in detail in \cite[Sect.\ VII.2]{Be68},
\cite{Fu76}, and a discussion of a class of canonical systems in
\cite[Ch.\ 8]{Sa97}.) In fact, at the conclusion of a meeting held in
honor of Professor Allan Krall at  the University of Tennessee-Knoxville
on October 10, 2002, Professor Krall noted in remarks, which he entitled
``Linear Hamiltonian systems involving difference equations'', that a
Weyl--Titchmarsh theory for general Hamiltonian systems of difference
equations has yet to be developed. By Hamiltonian  system of difference
equations he meant those systems that arise naturally as a
discretization of linear Hamiltonian systems of differential equations
(cf.\ \eqref{HSc}--\eqref{HSe}), and in analogy to the material developed
in \cite{HS81}--\cite{HS86}, the principal aim should be to construct the
matrix-valued Weyl--Titchmarsh function and develop the related
specral theory of such systems. In part, our paper is meant to follow up
on the challenge extended by Professor Krall and develop
Weyl--Titchmarsh theory for singular discrete Hamiltonian systems as a
natural extension of the existing theory for scalar Jacobi equations (cf.\
\cite{As66}, \cite[Sect.\ VII.1]{Be68}, \cite{GKT96}, \cite{GS97},
\cite{GT96}, \cite[Ch.\ 2]{Te00} and the references therein). The actual
model we follow closely in this paper is our recent treatment of
Dirac-type systems in \cite{CG02}.

In this paper we develop the basic theory of matrix-valued
Weyl--Titchmarsh M-functions and the associated Green's matrices for
whole-line and half-line self-adjoint Hamiltonian finite
difference systems defined as follows. Let $m\in\bbN$ and
\begin{equation}
B=\{B(k)\}_{k\in\bbZ}\subset\bbC(\bbZ)^{2m\times 2m}, \quad
\rho=\{\rho(k)\}_{k\in\bbZ}\subset\bbC(\bbZ)^{m\times m}
\end{equation}
with $\bbC(\bbZ)^{r\times s}$, the space of sequences of
complex $r\times s$ matrices, $r,s\in\bbN$, where $B(k)$ and $\rho(k)$ are
assumed to be self-adjoint and nonsingular matrices for all $k\in\bbZ$. We
denote by $S^\pm$ the shift operators acting upon
$\bbC(\bbZ)^{m\times s}$, that is,
\begin{equation}
S^{\pm}f(\,\cdot\,)=f^{\pm}(\,\cdot\,)=f (\, \cdot\,\pm 1), \quad
f\in \bbC(\bbZ)^{m\times s}.
\end{equation}
Moreover, let
\begin{equation}
A=\{A(k)\}_{k\in\bbZ}\subset\bbC(\bbZ)^{2m\times 2m},
\end{equation}
such that
\begin{equation}
A(k)=\begin{pmatrix}A_{1,1}(k)&A_{1,2}(k)
\\A_{2,1}(k)
& A_{2,2}(k) \end{pmatrix}\ge 0, \quad k\in\bbZ, \label{1.11}
\end{equation}
where
$A_{u,v}=\{A_{u,v}(k)\}_{k\in\bbZ}\in\bbC(\bbZ)^{m\times m}$, $u,v =1,2$.
Introducing the following linear difference expression
\begin{equation}\label{1.10}
\calS_{\rho} - B , \quad \calS_{\rho}= \begin{pmatrix}0 & \rho
S^+ \\ \rho^- S^- & 0 \end{pmatrix},
\end{equation}
the eigenvalue equation, or \emph{discrete Hamiltonian system} on the
whole-line considered in this paper, is then given by
\begin{equation}\label{1.12}
\calS_{\rho} \varPsi(z,k)=[zA(k)+B(k)]\varPsi(z,k), \quad
z\in\bbC, \; k\in\bbZ.
\end{equation}
Here $z$ plays the role of the spectral parameter and
\begin{equation}\label{1.13}
\varPsi(z,k) = \begin{pmatrix}\psi_1(z,k)\\
\psi_2(z,k) \end{pmatrix}, \quad \psi_j(z,\, \cdot\,)\in
\bbC(\bbZ)^{m\times r}, \; j=1,2
\end{equation}
with $1\leq r\leq 2m$, and $S^{\pm}\psi_j(z,\, \cdot\,) =
\psi_j(z,\,\cdot\,\pm  1)$, $j=1,2$. Analogously, we will consider
\eqref{1.12} on a half-line. Of course, at finite endpoints of the
underlying interval (and possibly also at the point(s) at infinity), the
formally self-adjoint Hamiltonian system \eqref{1.12} needs to be
supplied with appropriate self-adjoint boundary conditions to render it
self-adjoint. This will be discussed in Sections \ref{s2} and \ref{s3}. 

Forms such as \eqref{1.12} arise naturally when discretizing a Hamiltonian
system of first-order ordinary differential equations,
\begin{equation}\label{1.14}
J\varPsi'(z,x)=[z\calA(x) + \calB(x)]\varPsi(z,x), \quad x\in\bbR,
\end{equation}
as discussed in the next section (cf.\ the discussion
following \eqref{HSc}).

In Section \ref{s2} we set up the basic Weyl--Titchmarsh formalism
associated with \eqref{1.12}. We discuss possible normal forms of
\eqref{1.12} and show that $\rho$ can be assumed to be diagonal and positive
definite without loss of generality. Subsequently, we introduce the
necessary tools to discuss separated boundary conditions associated
with \eqref{1.12} on a finite interval and then define the corresponding
$m\times m$ matrix-valued Weyl--Titchmarsh function, the Weyl disk, and
the Weyl circle. The latter is shown to correspond to regular boundary
value problems associated with \eqref{1.12} on a finite interval with
separated self-adjoint boundary conditions at the endpoints. Next, the
Herglotz property of the Weyl--Titchmarsh function is established and
different boundary conditions at one endpoint (keeping the boundary
condition fixed at the other endpoint) are shown to be related by linear
fractional transformations. The typical nesting property of Weyl disks
associated with a finite interval then yield the existence of a limiting
Weyl disk as the finite interval approaches a half-line. The
limiting disk is nonempty, closed, and convex. The elements of the
limit disk turn out to be $m\times m$ matrix-valued Herglotz functions of
rank $m$. If the limiting Weyl disk consists of just a point, one then has
the important limit point case.

In our final Section \ref{s3} we then consider boundary value problems and
Green's functions associated with the discrete Hamiltonian system
\eqref{1.12} and appropriate self-adjoint boundary conditions on the
whole-line and on half-lines.

The results on Green's functions in Section \ref{s3} are fundamental for
the concrete applications we have in mind in our subsequent paper
\cite{CGR04}. There we will consider trace formulas and Borg-type
uniqueness theorems associated with matrix-valued Jacobi operators and
certain (supersymmetric) Dirac-type difference operators, which turn out
to be interesting special cases of the discrete Hamiltonian system
\eqref{1.12}. These special cases have interesting applications to
hierarchies of completely integrable nonabelian nonlinear evolution
equations. In fact, the matrix-valued Jacobi difference expression
\eqref{HSd} subject to \eqref{2.16} yields a Lax operator for the
nonabelian Toda hierarchy (cf., e.g., \cite[Sect.\ 12.2]{Te00},
\cite[Sects.\ 3.1, 3.2]{To89},) and the Dirac-type difference expression
\eqref{HSd} subject to \eqref{2.15} yields a Lax operator for the
nonabelian Kac--van Moerbeke hierarchy (cf., e,g., \cite{BGHT98},
\cite{GHSZ93}, \cite[Sect.\ 14.1]{Te00}, \cite[Sect.\ 3.8]{To89}).

\medskip
\noindent {\bf Dedication.} It is with great pleasure that we dedicate
this paper to Norrie Everitt on the occasion of his 80th birthday. His
enormous influence on the field of ordinary differential operators is
universally admired. In the very special context of this paper, we refer,
in particular, to his fundamental papers \cite{CE67},
\cite{Ev59}--\cite{EZ79}, which paved the way for a systematic treatment
of general Hamiltonian systems and inspired a whole generation of
scientists to enter this field. 

\section{Weyl--Titchmarsh Matrices for Finite Difference \\ Hamiltonian
Systems} \label{s2}

We now turn to the Weyl--Titchmarsh theory for Hamiltonian systems of
finite difference operators.
The model for this part of our discussion is the analogous development
of the theory presented in \cite{CG02} which in turn is based upon
the theory developed  by Hinton and Shaw in a series of papers devoted to 
the spectral theory of (singular) Hamiltonian systems of differential
equations \cite{HS81}--\cite{HS86} (see  also \cite{Kr89a}, \cite{Kr89b}).

Throughout this paper, matrices will be considered over the field
of complex numbers $\bbC$. With $M$ in the space of $r\times s$
complex matrices, $\bbC^{r\times s}$, $r,s\in\bbN$, let $M^\top$
denote the transpose,  and let $M^*$ denote the adjoint or
conjugate transpose of the matrix $M$. Let $M\ge 0$ and $M\le 0$
indicate that $M$ is nonnegative and nonpositive respectively.
Similarly, $M>0$ (respectively, $M<0$) denotes a positive definite
(respectively, negative definite) matrix. Moreover, let
$\Im(M)=(M-M^*)/(2i)$ and $\Re(M)=(M+M^*)/2$ denote the imaginary
and real parts  of the matrix $M$.

Denote by $\bbC(I)^{r\times s}$ the space of sequences, defined on
$I\subseteq\bbZ$, of complex $r\times s$ matrices where $1\le  s\le 2m$,
and where typically $r\in\{m, 2m\}$. Denote by $\ell^\infty(I)^{r\times
s}$ the sequence space of complex $r\times s$ matrices bounded on
$I\subseteq\bbZ$ with respect to the norm
$\|\cdot\|_{\ell^\infty(I)^{r\times s}}$, while $\ell^p(I)^{r\times s}$
denotes the space of sequences $p$-summable on $I\subseteq\bbZ$ with
respect to the norm $\|\cdot\|_{\ell^p (I)^{r\times s}}$. Let $S^\pm$
denote the shift operators on $\bbC(\bbZ)^{m\times r}$, that is,
\begin{equation}
S^{\pm   }f(\,\cdot\,) = f(\, \cdot\,\pm 1), \quad f^{\pm}= S^{\pm}f,
\quad f\in\bbC(\bbZ)^{m\times r}.
\end{equation}
Moreover, with $\Psi\in\bbC(\bbZ)^{2m\times r}$, let
\begin{equation}
\Psi= \begin{pmatrix}\psi_1\\  \psi_{2}
\end{pmatrix}, \quad
\hatt{\Psi}=\begin{pmatrix} \psi_1\\ \psi_2^+ \end{pmatrix}, \quad
\psi_j\in\bbC(\bbZ)^{m\times r},\quad j=1,2.
\end{equation}

Unless explicitly stated otherwise, $[c,d]\subset\bbR$ will mean the
discrete interval
$[c,d]\cap\bbZ$, with $c, d\in \bbZ$, possibly with $d<c$; the
trivial interval occurring when $c=d$. If $c\ne d$, let $^+[c,d]$
denote the discrete interval
\begin{equation}
^+[c,d]=[\min\{c,d\}+1,\max \{c,d\}].
\end{equation}
Evaluation may be expressed by
\begin{equation}
\psi {\big |}_c = \psi(c),
\end{equation}
while differences may be expressed by
\begin{equation}
\psi {\big |}_{[c,d]}= \psi(\max \{c,d\}) - \psi(\min\{c,d\}).
\end{equation}
Sums over discrete intervals may be expressed by
\begin{equation}
\sum_{k\in[c,d]}\psi(k) = \sum_{n=\min\{c,d\}}^{\max \{c,d\}}\psi(k).
\end{equation}
These conventions will turn out to be useful in connection with the
functional $E_\ell (M)$ introduced in \eqref{2.380} in the sense that they
permit us to avoid numerous case distinctions associated with $k_0>\ell$,
$k_0=\ell$, $k_0<\ell$, etc.

Next, let $m\in\bbN$, and let $J\in\bbC^{2m\times 2m}$, $J_\rho(k), 
I_\rho(k)\in\bbC(\bbZ)^{2m\times 2m}$ be defined for $k\in\bbZ$ by
\begin{equation}\label{d2.1a}
J=\begin{pmatrix}0& I_m \\ -I_m & 0  \end{pmatrix}, \hspace{5pt}
J_{\rho}(k)=\begin{pmatrix}0& \rho(k) \\ -\rho(k) & 0  \end{pmatrix},
\hspace{5pt}
I_{\rho}(k)=\begin{pmatrix}\rho(k) & 0 \\ 0 & \rho(k)  \end{pmatrix}.
\end{equation}
Here $I_m$ denotes the $m\times m$ identity matrix in $\bbC^m$ and
$\rho(k)\in\bbC^{m\times m}$ is self-adjoint and nonsingular for all
$k\in\bbZ$. Let $A_{u,v}(k),B_{u,v}(k)\in\bbC(\bbZ)^{m\times m}$ for $u,v
=1,2$ and $k\in\bbZ$. Moreover, for $k\in\bbZ$, let
\begin{align}\label{d2.1b}
A(k)&=\begin{pmatrix}A_{1,1}(k)&A_{1,2}(k)
\\A_{2,1}(k)
& A_{2,2}(k) \end{pmatrix}\ge 0, \\ \quad
B(k)
&=\begin{pmatrix}B_{1,1}(k)&B_{1,2}(k) \\B_{2,1}(k) &
B_{2,2}(k)
\end{pmatrix}=B(k)^* \label{d2.1c}.
\end{align}
In terms of the operator $\rho S^+$ and its formal adjoint $\rho^-
S^-$, let $\calS_{\rho}$ denote the formally  self-adjoint
matrix-valued difference expression given by
\begin{equation}\label{d2.1d}
\calS_{\rho}= \begin{pmatrix}0 & \rho S^+ \\ \rho^- S^- & 0 \end{pmatrix}.
\end{equation}

With $A(k)$, $B(k)$,  $\calS_{\rho}$ defined in \eqref{d2.1b} -- \eqref{d2.1d},
we consider the general difference expression given by
\begin{equation}\label{HSd}
\calS_{\rho} - B,
\end{equation}
and its associated eigenvalue equation, or  general \emph{Hamiltonian system},
given by
\begin{subequations}\label{HS}
\begin{equation}\label{HSa}
\calS_{\rho} \varPsi(z,k)=[zA(k)+B(k)]\varPsi(z,k), \quad
z\in\bbC, \; k\in\bbZ.
\end{equation}
Here $z$ plays the role of the spectral parameter and
\begin{equation}\label{HSb}
\varPsi(z,k) = \begin{pmatrix}\psi_1(z,k)\\
\psi_2(z,k) \end{pmatrix}, \quad \psi_j(z,\, \cdot\,)\in
\bbC(\bbZ)^{m\times r}, \; j=1,2
\end{equation}
with $1\leq r\leq 2m$, and $S^{\pm}\psi_j(z,\, \cdot\,) =
\psi_j(z,\,\cdot\,\pm  1)$, $j=1,2$.
Such a Hamiltonian system is said to be  \emph{well-posed} when it
possesses unique solutions defined for all $k\in\mathbb{ Z}$
associated with prescribed initial values of the type
\begin{equation}
\hatt\varPsi(z,k_0)\in\bbC^{2m}. \lb{HSf}
\end{equation}
\end{subequations}
A necessary and sufficient condition for well-posedness is given in
\eqref{2.17a} below.

For our discussion concerning the Weyl--Titchmarsh theory of
the Hamiltonian system \eqref{HSa}, we also adopt a definiteness
condition like that of Atkinson~\cite{At64}. We briefly sum up all
hypotheses on the coefficients in \eqref{HSa} as follows:

\begin{hypothesis}\label{h2.1}
We assume that our Hamiltonian system satisfies
\begin{align}
& A(k)\geq 0, \quad B(k)=B(k)^*, \quad \rho(k)>0 \, \text{ for all
$k\in\bbZ$}, \lb{2.17b} \\
&\text{$zA_{1,2}(k)+B_{1,2}(k)$ is invertible 
for all $k\in\bbZ$ and $z\in\bbC$,}  \label{2.17a}  
\end{align}
and for all nontrivial solutions $\varPsi\in\bbC^{2m}$ of
\eqref{HSa}, we suppose that
\begin{equation}\label{2.17}
\sum_{k\in[c,d]}\varPsi(z,k)^*A(k)\varPsi(z,k) > 0,
\end{equation}
for every nontrivial discrete interval $[c,d]\subset\bbZ$ in the
case of the whole-line (resp.,
$[c,d]\subset [k_0,\infty)$ or $[c,d]\subset (-\infty,k_0]$ for
some $k_0\in\bbZ$ in the case of half-lines).
\end{hypothesis}

\begin{remark} \lb{r2.2}
Of course, condition \eqref{2.17a} requires invertibility of
$B_{1,2}$ $($and hence that of $B_{2,1}$$)$. Moreover, it is
equivalent to  invertibility of
$zA_{2,1}(k)+B_{2,1}(k)$ for all $k\in\bbZ$ and $z\in\bbC$.
In addition, condition \eqref{2.17a} guarantees existence and
uniqueness of the initial value problem
\eqref{HSa},
\eqref{HSf} by an explicit step by step construction of the solution
$\varPsi(z,k)$,
$k\in\bbZ$: Given $\hatt\varPsi(z,k_0)$, one needs to invert
$[zA_{2,1}(z,\ell)+B_{2,1}(z,\ell)]$ for all $\ell\geq k_0+1$ to
construct $\varPsi(z,k)$ for all $k\geq k_0+1$ and one needs to invert
$[zA_{1,2}(z,\ell)+B_{1,2,}(z,\ell)]$ for all $\ell\leq k_0-1$ to
construct $\varPsi(z,k)$ for all $k\leq k_0$.
\end{remark}

To avoid numerous case distinctions we will suppose the whole-line part
of Hypothesis \ref{h2.1} throughout this section. We will make an
explicit distinction between the whole-line and half-line cases in Section
\ref{s3}. 

Forms such as \eqref{HSa} arise naturally when discretizing a Hamiltonian
system of first-order ordinary differential equations,
\begin{equation}\label{HSc}
J\varPsi'(z,x)=[z\calA(x) + \calB(x)]\varPsi(z,x), \quad x\in\bbR,
\end{equation}
by replacing $\varPsi'(z,x)$ with the difference
expression given by
\begin{equation}\
\left(\negthickspace \begin{array}{r}-\partial^* \psi_1(z,k) \\ \partial
\psi_2(z,k) \end{array}\negthickspace \right),
\end{equation}
where the formally adjoint operators $\partial$ and $\partial^*$ are defined by
\begin{equation}\label{HSe}
   \partial = S^+-I_m, \quad \partial^* = S^--I_m,
\end{equation}
and where $I_m$ represents the identity matrix in
$\bbC(\bbZ)^{m\times m}$. These forms also arise when considering
matrix-valued Jacobi
operators (cf.\ \cite{GKM02}), or when considering the matrix-valued
generalizations of the super-symmetric Dirac-type operators considered in
\cite{BGHT98},
\cite{GHSZ93}, and \cite[Sect.\ 14.1]{Te00}. In particular \eqref{HSd}
represents a super-symmetric Dirac-type operator in \eqref{d2.1c} when
$B_{11}(k)=B_{22}(k)=0$, that is,
\begin{equation}\label{2.15}
A(k)=I_{2m}, \quad B(k) = \begin{pmatrix} 0 & b(k)\\b(k)^* & 0
\end{pmatrix}, \quad k\in\bbZ,
\end{equation}
and is relevant to the Kac--van Moerbeke system.
Alternatively, \eqref{HSa} represents
\begin{equation}
\partial\, p\, \partial^* y+qy=zy, \label{2.15a}
\end{equation}
the matrix-valued Sturm--Liouville difference equation, when
\begin{equation}\label{2.16}
\rho(k)=I_{m}\ ,\quad A(k)=\begin{pmatrix}I_m &
0\\0&0\end{pmatrix}, \quad B(k)=\begin{pmatrix}-q(k) &
I_m \\I_m & p(k)^{-1} \end{pmatrix}, \quad k\in\bbZ,
\end{equation}
and $B(k)^*= B(k)$, $k\in\bbZ$.  Equation \eqref{2.15a} is intimately
related to the Jacobi operator $H$.
More precisely, introducing the matrix-valued Jacobi difference expression
$L$ by
\begin{equation}
L=aS^+ + a^- S^- + b,
\end{equation}
where
\begin{align}
a&=-p^+, \quad b=p^+ +p +q, \\
p&=-a^-, \quad q=a+a^- + b
\end{align}
are $m\times m$ matrices, the matrix-valued Sturm--Liouville difference
equation \eqref{2.15a} is equivalent to the equation
\begin{equation}
Ly=zy. \label{2.20}
\end{equation}

Note that the examples cited in \eqref{2.15} and in \eqref{2.16}
satisfy the requirements \eqref{2.17b} and \eqref{2.17} in Hypothesis
\ref{h2.1}. Moreover, \eqref{2.17a} is automatically satisfied for
the example described in \eqref{2.16}.

When considering the spectral or inverse spectral theory of
\eqref{HSd} (especially in the context of \cite{CGR04}), we may choose,
without loss of generality, a more restrictive normal
form of $\rho$ in which $\rho$ represents a diagonal and positive
definite matrix.

\begin{lemma} \lb{l2.3}
The difference expression in \eqref{HSd} is unitarily
equivalent to another such expression in which $\rho$ is diagonal
and positive definite.
\end{lemma}
\begin{proof}
Let $Q(k)\in\bbC^{m\times m}$ define a unitary matrix such that
$Q(k)\rho(k)Q(k)^{-1} = \wti d(k)$, where $\wti d(k)\in\bbR^{m\times
m}$ is diagonal and self-adjoint for all $k\in\bbZ$. Then,
\begin{equation}
{U_{\rho}}(\calS_{\rho} - B)U_{\rho}^{-1} = \calS_{\wti d}- \wti B,
\quad \wti B = {U_{\rho}}BU_{\rho}^{-1}, \quad
U_{\rho}=\begin{pmatrix}Q&0\\0&Q^-\end{pmatrix}.
\end{equation}
Next, let $\wti\epsilon (k)\in\bbR^{m\times m}$ be
a diagonal matrix for which
$(\wti\epsilon (k))_{\ell,\ell}\in\{+1,-1\}$,
$\ell=1,\dots,m$. Define $\epsilon(k)\in\bbR^{m\times m}$ by $\epsilon
(k)=\wti\epsilon (k) \wti \epsilon (k+1)$ and choose
$\wti\epsilon (k)$ so that $d=\epsilon \wti d >0$. Then,
\begin{equation}
U_{\epsilon }(\calS_{d} - \wti B)U_\epsilon^{-1}  =
\calS_{d} -\calB, \quad \calB = U_{\epsilon } \wti B
U_\epsilon^{-1}, \quad U_{\epsilon}
= \begin{pmatrix}\wti\epsilon &0\\0&\wti\epsilon \end{pmatrix},
\end{equation}
thus showing that $\calS_{\rho} - B$ is unitarily equivalent to a
difference expression of type \eqref{HSd} for which
$\rho$ is diagonal and positive definite.
\end{proof}

\begin{definition}\label{d2.4}
By a general difference  expression and its associated Hamiltonian system
of first-order difference equations, we mean \eqref{HSd} and \eqref{HS},
respectively,  subject to Hypothesis \ref{h2.1} as well as the
additional assumption that the matrix $\rho$ is positive definite.
\end{definition}

Thus, we assume the following set of assumptions for the remainder
of this paper:

\begin{hypothesis}\label{h2.4}
In addition to Hypothesis \ref{h2.1} assume that $\rho$ is 
positive definite.
\end{hypothesis}

Next, we introduce a set of matrices which will serve to describe
boundary data for
separated boundary conditions to be associated with the Hamiltonian system
given in \eqref{HSa}:

\begin{definition}\label{dBD}
Let $\mathcal{B}_d$ denote the set of matrices $\gamma = (\gamma_1\;
\gamma_2)$ with $\gamma_j \in \bbC^{m\times m}$,
$j=1,2$, which satisfy the following
conditions,
\begin{subequations}\label{BD}
\begin{equation}\label{BDa}
\rank (\gamma)  = m,
\end{equation}
and that either
\begin{equation}\label{BDc}
\Im (\gamma_1\gamma_2^*) \le 0,
\quad \text{or}
\quad \Im (\gamma_1\gamma_2^*) \ge 0,
\end{equation}
where $(2i)^{-1}\, \gamma J\gamma^*=\Im (\gamma_1\gamma_2^*)$.
Given the rank condition in \eqref{BDa},
we assume, without loss of generality in what follows, the
normalization
\begin{equation}\label{BDd}
\gamma\gamma^*  = I_m.
\end{equation}
With $\rho(k)\in\mathbb{C}^{m\times m}$ positive definite and diagonal,
and $I_{\rho}$ as given in
\eqref{d2.1a}, $\wti{\gamma}$ is given by
\begin{equation}\label{BDe}
\wti{\gamma}(k)=\gamma I_{\rho}(k)^{1/2}. 
\end{equation}
\end{subequations}
\end{definition}

In \eqref{BDe} (and in the remainder of this paper) $\rho^{1/2}$ will
always denote the unique positive definite square root of $\rho>0$.

\begin{remark} \label{r2.4}
With $\gamma\in\bbC^{m\times 2m}$, the conditions
\begin{equation}
\gamma\gamma^*=I_m, \quad \gamma J\gamma^*=0
\label{r2.4a}
\end{equation}
imply that $\gamma\in\calB_d$,
and
explicitly read
\begin{equation}
\gamma_1\gamma_1^*
+\gamma_2\gamma_2^*=I_m, \quad
\gamma_1\gamma_2^* -\gamma_2\gamma_1^*=0. \label{r2.4b}
\end{equation}
\end{remark}

In fact, from \eqref{r2.4a} one also obtains
\begin{equation}
\gamma_1^*\gamma_1 +\gamma_2^*\gamma_2=I_m, \quad
\gamma_1^*\gamma_2 -\gamma_2^*\gamma_1=0, \label{r2.4c}
\end{equation}
as is clear from
\begin{equation}
\begin{pmatrix} \gamma_1 & \gamma_2\\ -\gamma_2 & \gamma_1 \end{pmatrix}
\begin{pmatrix} \gamma_1^* & -\gamma_2^*\\ \gamma_2^* & \gamma_1^*
\end{pmatrix}=I_{2m}=\begin{pmatrix} \gamma_1^* & -\gamma_2^*\\
\gamma_2^* & \gamma_1^* \end{pmatrix}\begin{pmatrix} \gamma_1 & \gamma_2\\
-\gamma_2 & \gamma_1 \end{pmatrix}, \label{r2.4d}
\end{equation}
since any left inverse matrix is also a right inverse, and vice
versa. Moreover, from \eqref{r2.4c} or \eqref{r2.4d}, we obtain
\begin{equation}\label{r2.4e}
\gamma^*\gamma J + J\gamma^*\gamma =  J.
\end{equation}

With $\alpha\in\bbC^{m\times 2m}$ satisfying \eqref{r2.4a} and
with $\ta=\ta(k_0)$ defined according to \eqref{BDe}, let
$\Psi(z,k,k_0,\ta)$ denote a normalized fundamental system of
solutions for the Hamiltonian system  \eqref{HSa} described in
Definition \ref{d2.4} which for some $k_0\in\bbZ$ satisfies
\begin{subequations}\label{FS}
\begin{equation}\label{FSa}
\hatt{\Psi}(z,k_0,k_0,\wti\alpha)=I_{\rho}(k_0)^{-1}(\wti\alpha^* \
J\wti\alpha^*)=I_{\rho}(k_0)^{-1/2}(\alpha^*\ J\alpha^*)  .
\end{equation}
We partition $\Psi(z,k,k_0,\wti\alpha)$ as follows,
\begin{align}
\Psi(z,k,k_0,\wti\alpha)&=(\Theta(z,k,k_0,\wti\alpha)\;\;\;
\Phi(z,k,k_0,\wti\alpha))\label{FSb}\\
&=\begin{pmatrix}\theta_1(z,k,k_0,\wti\alpha) &
\phi_1(z,k,k_0,\wti\alpha)\\
\theta_2(z,k,k_0,\wti\alpha)& \phi_2(z,k,k_0,\wti\alpha)
\end{pmatrix}, \label{FSc}
\end{align}
\end{subequations}
where $\theta_j(z,k,k_0,\wti\alpha)$ and $\phi_j(z,k,k_0,\wti\alpha)$
for $j=1,2$ are $m\times m$ matrices, entire with
respect to $z\in\bbC$, and normalized according to
\eqref{FSa}. One can now prove the following result.

\begin{lemma}\label{l2.4}
Let $\Theta(z,k)=\Theta(z,k,k_0,\wti\alpha)$ and
$\Phi(z,k)=\Phi(z,k,k_0,\wti\alpha)$ be defined in \eqref{FS} with
$\alpha,\beta\in\calB_d$, and with
$\Im(\alpha_1\alpha_2^*)=0$. Let $\wti\alpha=\wti\alpha(k_0)$ and
$\wti\beta=\wti\beta(\ell)$ be defined according to \eqref{BDe}. Then, for
$\ell\ne k_0$, $\wti\beta\hatt{\Phi}(z,\ell)$ is singular  if and only
if $z$ is an eigenvalue for the  regular boundary value problem
given by \eqref{HSa} together with the separated boundary
conditions
\begin{equation}\label{BC}
\wti\alpha\hatt{\varPsi}(z,k_0)=0, \quad \wti\beta\hatt{\varPsi}(z,\ell)=0.
\end{equation}
\end{lemma}

\noindent One observes that both regular boundary conditions described in
\eqref{BC} are self-adjoint when $\Im(\beta_1\beta_2^*)=0$.

In light of Lemma \ref{l2.4}, it is possible to introduce, under
appropriate conditions, the $m\times m$ matrix-valued meromorphic
function, $M\big(z,\ell,k_0,\wti\alpha,\wti\beta\big)$, as follows.

\begin{definition}\label{dMF}
Let \eqref{FS} define $ \Theta(z,k,k_0,\wti\alpha)$,
and $\Phi(z,k,k_0,\wti\alpha)$ with
$\alpha,\beta\in\calB_d$, and $\Im(\alpha_2\alpha_1^*)=0$.  For $\ell\ne
k_0$, and $\wti\beta\hatt{\Phi}(z,\ell,k_0,\wti\alpha)$ nonsingular,
define
\begin{equation}\label{MF}
M\big(z,\ell,k_0,\wti\alpha,\wti\beta\big) =
-[\wti\beta\hatt{\Phi}(z,\ell,k_0,\wti\alpha)]^{-1}
[\wti\beta\hatt{\Theta}(z,\ell,k_0,\wti\alpha)].
\end{equation}
$M\big(z,\ell,k_0,\wti\alpha,\wti\beta\big)$ is said to be the
{\it Weyl--Titchmarsh $M$-function} for the regular boundary value
problem described in Lemma \ref{l2.4}.
\end{definition}

By means of the equations
\begin{equation}
v_1(z,k) = \rho(k)^{1/2}\psi_1(z,k), \quad v_2(z,k) =
\rho^-(k)^{1/2}\psi_2(z,k), \label{2.28a}
\end{equation}
and
\begin{equation}
\wti\alpha = \wti\alpha (k_0)=\alpha I_{\rho}(k_0)^{1/2}, \quad
\wti\beta = \wti\beta (\ell)=\beta I_{\rho}(\ell)^{1/2}, \label{2.36}
\end{equation}
the boundary value problem described in Lemma \ref{l2.4} is transformed
into one described by
\begin{equation}\label{2.300}
\calS_{I_m} V(z,k)=[z\wti{A}(k)+\wti{B}(k)]V(z,k), \quad
V(z,k)=\begin{pmatrix} v_1(z,k) \\ v_2(z,k) \end{pmatrix},
\end{equation}
\begin{equation}
\alpha\hatt{V}(z,k_0)=0, \quad \beta\hatt{V}(z,\ell)=0,
\end{equation}
where
\begin{equation}
\wti{A}= \mathcal{D}^{-1/2}A\mathcal{D}^{-1/2}, \quad
\wti{B}= \mathcal{D}^{-1/2}B\mathcal{D}^{-1/2}, \quad
\mathcal{D}=\begin{pmatrix} \rho & 0\\0 & \rho^- \end{pmatrix}.
\end{equation}

The following remark will play a role in connection with inverse
spectral theory considerations in \cite{CGR04}.

\begin{remark}
$(i)$ In general, \eqref{2.28a} and \eqref{2.36} do not define a unitary
transformation. However, if $\Psi_{I_m}(z,k,k_0,\alpha)$ represents the
fundamental solution matrix of \eqref{2.300} such that
$\Psi_{I_m}(z,k_0,k_0,\alpha)= (\alpha^* \, J\alpha^*)$ and if
\begin{equation}
M_{I_m}(z,\ell,k_0,\alpha,\beta)=
-[\beta\hatt{\Phi}_{I_m}(z,\ell,k_0,\alpha)]^{-1}
[\beta\hatt{\Theta}_{I_m}(z,\ell,k_0,\alpha)]
\end{equation}
represents the associated
M-function for \eqref{2.300}, while
$M\big(z,\ell,k_0,\wti\alpha,\wti\beta\big)$ is defined in \eqref{MF},
then
\begin{equation}
M_{I_m}(z,\ell,k_0,\alpha,\beta)
=M\big(z,\ell,k_0,\wti\alpha,\wti\beta\big).
\end{equation}
So one can replace $\rho$ by $I_m$ and $\wti\alpha,\,
\wti\beta$ by $\alpha,\, \beta $ but possibly at the expense of more
complex expressions for $A$ and $B$. \\
$(ii)$ Assume that $M\big(z,\ell,k_0,\wti\alpha,\wti\beta\big)$ is
the $M$-function defined in \eqref{MF} $($with $\rho(k)>0$,
$k\in\bbZ$$)$ and $M_{\wti d}(z,\ell,k_0,\gamma,\delta)$ is the
$M$-function corresponding to 
\begin{align}
&\calS_{\wti d}\varPhi(z,k)=[z\wti A(k)+ \wti B(k)]\varPhi(z,k), \\
& \wti A(k) = U_{\rho}(k)A(k)U_{\rho}(k)^{-1}, \quad \wti B(k) =
{U_{\rho}(k)}B(k)U_{\rho}(k)^{-1}, \\
& U_{\rho}(k)=\begin{pmatrix}Q(k)&0\\0&Q(k)^-\end{pmatrix}, \quad
k\in\bbZ.
\end{align}
with $\wti d(k)>0$ a diagonal matrix as discussed in the proof of
Lemma \ref{l2.3} $($and $Q(k)$ a unitary $m\times m$ matrix as in the
proof of Lemma \ref{l2.3}, except that $\rho(k)$ was not assumed to be
positive definite in Lemma \ref{l2.3}$)$. Then,
\begin{equation}
M_{\wti d}(z,\ell,k_0,\gamma,\delta)
=M\big(z,\ell,k_0,\wti\alpha,\wti\beta\big),
\end{equation}
where 
\begin{equation}
\gamma=\wti\alpha\begin{pmatrix} Q(k_0)^{-1} & 0 \\ 0 & Q(k_0)^{-1}
\end{pmatrix},
\quad 
\delta=\wti\beta\begin{pmatrix} Q(\ell)^{-1} & 0 \\ 0 & Q(\ell)^{-1}
\end{pmatrix}.
\end{equation}
\end{remark}

\medskip

The Weyl--Titchmarsh $M$-function in \eqref{MF} is an $m\times m$
matrix-valued function with meromorphic entries whose poles
correspond to eigenvalues for the regular boundary value problem
given by the difference equation \eqref{HSa} and the
boundary conditions \eqref{BC}. Moreover, given the normalized
fundamental matrix, $\Psi(z,k,k_0,\ta)$, defined in \eqref{FS} and given
$M\in\bbC^{m\times m}$, one defines
\begin{equation}\label{2.19}
U(z,k,k_0,\ta)= \begin{pmatrix}
u_1(z,k,k_0,\ta)\\u_2(z,k,k_0,\ta) \end{pmatrix}=
\Psi(z,k,k_0,\ta)\begin{pmatrix}I_m\\M\end{pmatrix},
\end{equation}
with $u_j(z,k,k_0,\ta)\in \bbC^{m\times m}$, $j=1,2$. Then
$U(z,k,k_0,\ta)$ will satisfy the boundary condition at $k=\ell$ in
\eqref{BC} when $M=M\big(z,\ell,k_0,\ta,\tb\big)$. Intimately connected
with the matrices introduced in Definition \ref{dMF} is the set of
$m\times m$ complex matrices known as the Weyl disk. Several
characterizations of this set have appeared in the literature for
the Hamiltonian system of differential equations given in
\eqref{HSc} (see, e.g., \cite{At64}--\cite{At88a},
\cite{HS93}, \cite{HS81}, \cite{Kr89a}, \cite{Or76}). By analogy,
such definitions also exist for the Hamiltonian difference equation
\eqref{HSa}.

To describe this set, we first introduce the matrix-valued function
$E_\ell (M)$: With $\ell\ne k_0$, $z\in\bbC\backslash\bbR$, and with
$U(z,\ell,k_0,\ta)$  defined by \eqref{2.19} in terms of a  matrix
$M\in\bbC^{m\times m}$, let
\begin{equation}\label{2.380}
E_\ell(M) = \sigma(\ell,k_0,z)\hatt{U}(z,\ell,k_0,\ta)^*
(iJ_{\rho}(\ell)) \hatt{U}(z,\ell,k_0,\ta),
\end{equation}
where
\begin{equation}
\sigma(\ell,k,z)=\frac{(\ell-k)\Im (z)}{|(\ell-k)\Im (z)|},\quad
\sigma(\ell,k)=\sigma(\ell,k,i)  
\end{equation}
with $\ell\ne k$, and  $\ell,k\in\bbZ$.

\begin{definition}\label{dWD}
Let the following be fixed: Integers $k_0$ and $\ell\ne k_0$, 
   $\alpha\in\calB_d$, and $z\in\bbC\backslash\bbR$.
$\calD(z,\ell,k_0,\ta)$ will denote the collection of all $M\in
\bbC^{m\times m}$ for which $E_\ell (M)\le 0$, where $E_\ell (M)$
is defined in \eqref{2.380}. $\calD(z,\ell,k_0,\ta)$ is said to be
a {\it Weyl disk}.   The set of $M\in \bbC^{m\times m}$ for which
$E_\ell (M) = 0$ is said to be a {\it Weyl circle} $($even when
$m>1$$)$. The interior of the Weyl disk is denoted by
$\calD(z,\ell,k_0,\ta)^{\circ}$.
\end{definition}

This definition leads to a representation that is a generalization
of the description first given by Weyl \cite{We10} in the context
of Sturm-Liouville differential expressions: a representation in
which $\calD(z,\ell,k_0,\ta)$ is homeomorphic to the set of
contractive matrices, that is, those matrices $V\in\bbC^{m\times m}$
for which $VV^*\le I_m$. This provides the justification for the
geometric terms of circle and disk (cf., e.g., \cite{HS93}, \cite{HS81},
\cite{Kr89a}, \cite{Or76}). From this representation
it is also seen that the \emph{interior} of the Weyl disk is nonempty
and corresponds to the collection of
all $M\in\bbC^{m\times m}$ for which $E_\ell (M)< 0$.

We next discuss some basic properties associated with elements of
the disk $\calD(z,\ell,k_0,\ta)$.  To this end, we introduce the
assumptions contained in the
next hypothesis for the parameters $k_0$ and $\ell$:

\begin{hypothesis}\label{h2.6}
If for the Hamiltonian system satisfying Hypothesis \ref{h2.4} it is given
that $A(k)>0$ for $k\in \mathbb{Z}$, then $k_0\ne \ell$; otherwise,
we assume that $^+[k_0, \ell]$ is nontrivial.
\end{hypothesis}
In the next lemma,
we note that the Weyl circle corresponds to the regular boundary
value problems  with separated, self-adjoint
boundary conditions described in in Lemma \ref{l2.4}.
This lemma is the analog in our discrete
setting of  Lemma 2.8 in \cite{CG02}. For convenience of the reader,
and to achieve a reasonable level of completeness, we produce the
corresponding short proof below.
\begin{lemma}\label{l2.11}
Given Hypothesis \ref{h2.6}, let $M\in\bbC^{m\times m}$,
and let $z\in\bbC\backslash\bbR$. Then, $E_\ell (M)=0$ if and only if
there is a $\beta\in \bbC^{m\times 2m}$ satisfying \eqref{r2.4a}
such that
\begin{equation}\label{2.43}
0=\tb\hatt{U}(z,\ell,k_0,\ta),
\end{equation}
where $U(z,\ell,k_0,\ta)$ is defined in \eqref{2.19} in
terms of $M$,  and where $\tb=\tb(\ell)$.  With $\beta$ so defined,
\begin{equation}\label{2.44}
M=-[\tb\hatt{\Phi}(z,\ell,k_0,\ta)]^{-1}
[\tb\hatt{\Theta}(z,\ell,k_0,\ta)],
\end{equation}
that is, $M=M\big(z,\ell,k_0,\ta,\tb\big)$. Moreover, $\beta\in\calB_d$
and may be chosen to satisfy \eqref{BDd}.
\end{lemma}
\begin{proof}
Let $z\in\bbC\backslash\bbR$, and suppose for a given
$M\in\bbC^{m\times m}$ that there is a $\beta\in\bbC^{m\times 2m}$
which satisfies \eqref{r2.4a} and such that \eqref{2.43} is
satisfied. Given that $\beta J \beta^* =0$ and that $\rank
(\tb)=\rank (I_{\rho}(\ell)^{-1}J\tb^*)=m$, then by \eqref{2.43}
there is a nonsingular $C\in \bbC^{m\times m}$ such
that
$\hatt{U}(z,\ell,k_0,\ta) = I_{\rho}(\ell)^{-1}J\tb^*C$.
Hence,
$E_\ell (M)=-i\sigma(\ell,k_0,z)C^*\beta J\beta^*C=0$.

Upon showing that $\tb \hatt{\Phi}(z,\ell,k_0,\ta)$ is nonsingular,
\eqref{2.44} will then follow from \eqref{2.43}. If
$\hatt{\Phi}(z,\ell)=\hatt{\Phi}(z,\ell,k_0,\ta)$ and
$\tb\hatt{\Phi}(z,\ell)$ is singular, then there are nonzero vectors
$v, w \in \bbC^{m}$ such that $\tb\hatt{\Phi}(z,\ell)v=0$, and such
that $\hatt{\Phi}(z,\ell)v = I_{\rho}(\ell)^{-1}J\tb ^*w$. Let
$\varPsi_j=\varPsi_j(z,k)$, $j=1,2$, denote solutions of
\eqref{HSa} with $z=z_j$, $j=1,2$. Noting that
\begin{equation}\label{2.21}
\partial (\hatt{\varPsi}_1^* J_{\rho}\hatt{\varPsi}_2)^- =
\varPsi_1^*\mathcal{S}_{\rho}\varPsi_2
-\left( \mathcal{S}_{\rho}\varPsi_1 \right)^*\varPsi_2
= (z_2 - \bar{z}_1) \varPsi_1^* A \varPsi_2.
\end{equation}
and recalling that $\Phi(z,\,\cdot\,)=\Phi(z,\,\cdot\,,k_0,\ta)$
is defined in \eqref{FS},  we  obtain
\begin{equation}
\hatt{\Phi}^* J_{\rho}\hatt{\Phi}{\big |}_{[k_0,\ell]}= (z -
\bar{z})\sum_{^+[k_0,\ell]} \Phi^* A \Phi.
\end{equation}
   by \eqref{2.21}.  Since $\hatt{\Phi}(z,k_0)^*
J_{\rho}(k_0)\hatt{\Phi}(z,k_0)=0$, we then see that
\begin{align}\label{2.24}
2i\Im (z)\sum_{k\in^+[k_0,\ell]} v^*\Phi^*(z,k) A(k)\Phi(z,k)v &=
\sigma(\ell,k_0)v^*\hatt{\Phi}(z,\ell)^*
J_{\rho}(\ell)\hatt{\Phi}(z,\ell)v \\
&=w^*\beta J \beta^*w =0.
\end{align}
By Hypothesis \ref{h2.6}, $\Im(z)=0$. This contradicts the
assumption that $z\in\bbC\backslash\bbR$.

Conversely, if $E_\ell (M)=0$ for a given $M\in \bbC^{m\times m}$,
then for $z\in\bbC\backslash\bbR$  let $\tb =[I_m\;
M^*]\hatt{\Psi}(z,\ell,k_0,\ta)^*J_{\rho}(\ell) =
\hatt{U}(z,\ell,k_0,\ta)^* J_{\rho}(\ell)$, and let $\beta=\tb
I_{\rho}(\ell)^{-1/2}$. Thus \eqref{2.43} is satisfied and   $\rank
(\beta) =\rank (\tb)=m$. Moreover,
$0=E_\ell (M)=2\sigma(\ell,k_0,z)\Im(\beta_1\beta_2^*)$. If for this
choice of $\beta$,  \eqref{BDd} is not yet satisfied  Let $\delta
= (\beta\beta^*)^{-1/2}\beta$, and $\wti{\delta}= \delta
I_{\rho}(\ell)^{1/2}$. Note that $0=\wti{\delta} U(z,\ell,k_0,\ta)$,
that $\Im(\delta_1\delta_2^*)= (\beta\beta^*)^{-1/2}
\Im(\beta_1\beta_2^*) (\beta\beta^*)^{-1/2}$, and hence that
$\delta\in\calB_d$.
\end{proof}

Next, we observe that a fundamental property holds for  matrices
in $\calD(z,\ell,k_0,\ta)$.

\begin{lemma} \label{l2.8}
Given Hypothesis \ref{h2.6}, let $M\in\bbC^{m\times m}$, and let
$z\in\bbC\backslash\bbR$. Then,
\begin{equation}\label{Hgz}
\sigma(\ell,k_0,z) \Im (M)>0,
\end{equation}
whenever $M\in\calD(z,\ell,k_0,\ta)$.  Moreover, whenever
$\beta\in\bbC^{m\times 2m}$ satisfies \eqref{r2.4a},
\begin{equation}\label{2.270}
M\big(\bar z,\ell,k_0,\ta,\tb\big)= M\big(z,\ell,k_0,\ta,\tb\big)^*.
\end{equation}
\end{lemma}
\begin{proof}
By \eqref{2.21},
\begin{subequations}\label{2.32}
\begin{align}
2i\Im(z)\sum_{^+[k_0,\ell]} U^* A U &=
\hatt{U}^* J_{\rho}\hatt{U}{\big |}_{[k_0,\ell]}\label{2.32a}\\
&=2i\sigma(\ell,k_0)\Im(M) + \sigma(\ell,k_0)\hatt{U}^* J_{\rho}\hatt{U}
{\big |}_\ell
\end{align}
\end{subequations}
with $U=U(z,\,\cdot\,,k_0,\ta)$ defined in \eqref{2.19}. Moreover, by
the definition of $E_\ell (M)$ in \eqref{2.380}, one obtains
\begin{equation}\label{2.390}
2\sigma(\ell,k_0,z)\Im(M) = -E_\ell (M)
+ 2|\Im(z)|\sum_{^+[k_0,\ell]} U^*AU.
\end{equation}
By Hypothesis \ref{h2.6} and Definition \ref{dWD}, one infers that
$\sigma(\ell,k_0,z) \Im (M)>0$.

To prove \eqref{2.270}, we first let
$\Psi(z)=\Psi(z,\,\cdot\,,k_0,\ta)$, where $\Psi$ is defined in
\eqref{FS}. By \eqref{2.21} we note that
$\hatt{\Psi}(\bar{z})^*J_{\rho}\hatt{\Psi}(z)=-J$. As a
consequence,
\begin{equation}
I_{\rho}^{1/2}J\hatt{\Psi}(z)(\hatt{\Psi}(\bar{z})J)^*I_{\rho}^{1/2}
=-I_{2m},
\end{equation}
and hence
$\hatt{\Psi}(z)J\hatt{\Psi}(\bar{z})^*=-JI_{\rho}^{-1}$.
We obtain
\begin{subequations}
\begin{align}
(\tb\hatt{\Phi}(z))(\tb\hatt{\Theta}(\bar{z}))^*{\big |}_\ell -
(\tb\hatt{\Theta}(z))(\tb\hatt{\Phi}(\bar{z}))^*{\big |}_\ell &=
-\tb JI_{\rho}^{-1}\tb^*{\big |}_\ell \\
&= -\beta J\beta^*=0.
\end{align}
\end{subequations}
Equation~\eqref{2.270} then follows immediately.
\end{proof}

For $\ell>k_0$, the function $M\big(z,\ell,k_0,\ta,\tb\big)$, defined by
\eqref{MF}, whose values satisfy \eqref{Hgz}, thus represents a
matrix-valued {\it Herglotz} function of rank $m$. Hence, for
$\Im(\beta_2\beta_1^*)=0$, poles of $M\big(z,\ell,k_0,\ta,\tb\big)$,
$\ell>k_0$, are at most of first-order, are real, and have
nonpositive residues. For $\ell<k_0$, $-M\big(z,\ell,k_0,\ta,\tb\big)$ is
a Herglotz matrix. Thus, one obtains a representation of
$\sigma(\ell,k_0)M\big(z,\ell,k_0,\ta,\tb\big)$ of the form (cf.\ 
\cite{Ca76}, \cite{GT00}, \cite{HS81}, \cite{HS82}, \cite{HS86})
\begin{align}\label{NP}
\sigma(\ell,k_0)M\big(z,\ell,k_0,\ta,\tb\big)&=
C_1\big(\ell,k_0,\ta,\tb\big) + zC_2\big(\ell,k_0,\ta,\tb\big)\no\\
&\quad +\int_{-\infty}^\infty
d\Omega\big(\lambda,\ell,k_0,\ta,\tb\big)\,\left(
\frac{1}{\lambda-z} -\frac{\lambda}{1+\lambda^2} \right), 
\end{align}
where $C_2\big(\ell,k_0,\ta,\tb\big)\ge 0$ and 
$C_1\big(\ell,k_0,\ta,\tb\big)$ are self-adjoint $m\times m$ matrices,
and where
$\Omega\big(\lambda,\ell,k_0,\ta,\tb\big)$ is a nondecreasing $m\times m$
matrix-valued function such that
\begin{subequations}\label{Mrep}
\begin{align}
&\int_{-\infty}^{\infty}\|
d\Omega\big(\lambda,\ell,k_0,\ta,\tb\big)\|_{\bbC^{m\times m}}\,
(1+\lambda^2)^{-1}  < \infty, \label{NPa} \\
&\Omega\big((\lambda, \mu],\ell,k_0,\ta,\tb \big)=
 \lim_{\delta\downarrow 0}
\lim_{\epsilon \downarrow 0}\frac{1}{\pi}\int_{\lambda+ \delta}^{\mu
+ \delta }d\nu\,  \Im\big[\sigma(\ell,k_0) M\big(\nu
+i\epsilon,\ell,k_0,\ta,\tb\big)\big]. \label{NPb}
\end{align}
\end{subequations}

In general, for self-adjoint boundary value problems,
$\Omega\big(\lambda,\ell,k_0,\ta,\tb\big)$ is piecewise constant with jump
discontinuities precisely at the eigenvalues of the boundary value
problem, and that in the matrix-valued Schr\"odinger and
Dirac-type cases $C_2=0$ in \eqref{NP}. Analogous statements apply
to $-M\big(z,\ell,k_0,\ta,\tb\big)$ if $\ell<k_0$. For such problems, we
note in the subsequent lemma that for fixed $\beta$, varying the
boundary data $\alpha$ produces Weyl--Titchmarsh matrices
$M\big(z,\ell,k_0,\ta,\tb\big)$ which are related  to each other by a
linear fractional transformations (see also \cite{GMT98},
\cite{GT00} for a general approach to such linear fractional
transformations).

\begin{lemma}\label{l2.9} 
Assume Hyothesis \ref{h2.4}. 
If $\alpha, \beta, \gamma\in\bbC^{m\times 2m}$ satisfy
\eqref{r2.4a} and if $\ta=\ta(k_0)$, $\tg=\tg(k_0)$, $\tb=\tb(\ell)$,
then
\begin{equation}
M_{\ta}= [-\alpha J \gamma^* + \alpha\gamma^* M_{\tg}]
[\alpha\gamma^* +  \alpha J\gamma^*M_{\tg}]^{-1}, \label{2.360}
\end{equation}
where $M_{\ta}=M\big(z,\ell,k_0,\ta,\tb\big)$, and
$M_{\tg}=M\big(z,\ell,k_0,\tg ,\tb\big)$.
\end{lemma}

\begin{remark}
{}From the proof of the previous lemma one infers, in general,  that
\begin{equation}
U_{\tg}(z,k) =
U_{\ta}(z,k)(\alpha \gamma^* + \alpha J \gamma^* M_{\tg} ).
\end{equation}
Moreover,  if $\alpha_0 =(I_m\; 0)$ and $\gamma_0=(0\ I_m)$ one
observes, in particular, that
\begin{equation}
M\big(z,\ell,k_0,\ta_0,\tb\big)=-M\big(z,\ell,k_0,\tg_0,\tb\big)^{-1}.
\end{equation}
\end{remark}

We  further note that the sets $\calD(z,\ell,k_0,\ta)$ are closed,
and convex, (cf., e.g., \cite{HS93}, \cite{HS84}, \cite{Kr89a},
\cite{Or76}). Moreover, by \eqref{2.390} and
Hypothesis \ref{h2.4}, one concludes that $E_\ell (M)$ is increasing.
This fact  implies that, as a function of $\ell$, the sets
$\calD(z,\ell,k_0,\ta)$ are  nesting in the sense that
\begin{equation}\label{2.28}
\calD(z,\ell_2,k_0,\ta)\subset \calD(z,\ell_1,k_0,\ta) \quad
\text{for}\quad k_0<\ell_1< \ell_2\quad \text{or} \quad \ell_2< \ell_1<k_0.
\end{equation}
Hence, the intersection of this nested sequence, as $\ell\to
\pm \infty$, is nonempty, closed and convex.  We say that this
intersection is a limiting set for the nested sequence.

\begin{definition}\label{dLWD}
Let  $\calD_\pm(z,k_0,\ta)$ denote the closed, convex set in the space
of $m\times m$ matrices which is the limit, as
$\ell\to \pm\infty$,  of the nested collection of sets
$\calD(z,\ell,k_0,\ta)$ given in Definition \ref{dWD}.
$\calD_\pm(z,k_0,\ta)$ is said to be a limiting {\em disk}.
Elements of $\calD_\pm(z,k_0,\ta)$ are denoted by $M_\pm(z,k_0,\ta)\in
\bbC^{m\times m}$.
\end{definition}

In light of the containment described in \eqref{2.28} and
Hypothesis \ref{h2.4}, for $\ell\ne k_0$ and $z\in\bbC\backslash\bbR$,
\begin{equation}\label{2.321}
\calD_\pm(z,k_0,\ta)\subset \calD(z,\ell,k_0,\ta),
\end{equation}
with emphasis on strict containment of the disks in \eqref{2.321}.
Moreover, by \eqref{2.390},
\begin{equation}\label{2.320}
M\in \calD_\pm(z,k_0,\ta) \text{ precisely when }E_\ell (M)<0
\text{ for all }  \ell \in(k_0, \pm\infty).
\end{equation}
In the next lemma, the interior points of the Weyl
disk are characterized in terms of certain elements of $\calB_d$ (cf.\
\eqref{BDe}). This lemma is the analog in our discrete setting
both in its statement and proof of Lemma 2.13 of \cite{CG02}.
\begin{lemma}\label{l2.12}
Given Hypothesis \ref{h2.6}, let $M\in\bbC^{m\times m}$,
and let $z\in\bbC\backslash\bbR$. Then, $E_\ell (M)<0$ if and only if
there is a $\beta\in \bbC^{m\times 2m}$   satisfying the condition
\begin{equation}\label{2.27a}
\sigma(\ell,k_0,z) \Im(\beta_1\beta_2^*)>0,
\end{equation}
and such that \eqref{2.43} holds
with $u_j(z,\ell)= u_j(z,\ell,k_0,\ta)$, $j=1,2$, defined in
\eqref{2.19} in terms of $M$. With $\beta$ so defined,
\eqref{2.44} holds, that is, $M=M\big(z,\ell,k_0,\ta,\tb\big)$. Moreover,
$\beta\in\calB_d$ and $\beta$ may be chosen to satisfy \eqref{BDd}.
\end{lemma}
\begin{proof}
Let $z\in \bbC\backslash\bbR$, and for a
given $M\in \bbC^{m\times m}$  suppose that there is a
$\beta\in\bbC^{m\times 2m}$ satisfying \eqref{2.27a} such that
\eqref{2.43} holds. The matrices $\beta_j$, $j=1,2$, are
invertible by \eqref{2.27a}, and by \eqref{2.43} it follows that
\begin{equation}\label{2.46}
\hatt{U}(z,\ell)=\begin{pmatrix}
-\tb_1^{-1}\tb_2\\I_m \end{pmatrix}
u_2^+(z,\ell,k_0,\ta).
\end{equation}
By \eqref{2.380} and \eqref{2.46}, we see that
\begin{equation}\label{2.47}
E_\ell (M) = -2\sigma(\ell,k_0,z) u_2^+(z,\ell)^*\rho^{1/2}\beta_1^{-1}
\Im (\beta_1\beta_2^*){(\beta_1^*)}^{-1}\rho^{1/2} u_2^+(z,\ell),
\end{equation}
and hence that $E_\ell (M)<0$ whenever \eqref{2.27a} holds.

Upon showing that $\tb\hatt{\Phi}(z,\ell)$ is nonsingular,
\eqref{2.44} will follow from \eqref{2.43}. If $\tb\hatt{\Phi}(z,\ell)$
is singular, then there is a nonzero vector $v\in \bbC^{m}$ such that
$\tb\hatt{\Phi}(z,\ell)v=0$. By the nonsingularity of $\beta_j$, $j=1,2$,
$\phi_1(z,\ell)v = -\tb_1^{-1}\tb_2\phi_2^+(z,\ell)v$, and
as a result, \eqref{2.24} yields
\begin{align}
&2|\Im (z)|\sum_{k\in^+[k_0,\ell]} v^*\Phi^*(z,k) A(k)\Phi(z,k)v \no \\
&=-2\sigma(\ell,k_0,z) v^*u_2^+(z,\ell)^*\rho^{1/2}\beta_1^{-1}
\Im (\beta_1\beta_2^*){(\beta_1^*)}^{-1}\rho^{1/2} u_2^+(z,\ell)v,
\end{align}
and hence a contradiction given \eqref{2.27a} (cf. \eqref{2.17}).

Conversely, if $E_\ell (M)=2\sigma(\ell,k_0,z)\Im(u_1^*\rho
u_2^+)\big|_\ell <0$ for a given $M\in\bbC^{m\times m}$, then for
$z\in\bbC\backslash\bbR$,  $u_2^+(z,\ell,k_0,\ta)$ is nonsingular.
Let $\beta_1 =I_{\rho}(\ell)^{-1/2}$ and let
$\beta_2=-u_1({u_2}^+)^{-1}I_{\rho}^{-1/2}\big|_\ell$. Then, for
this choice of $\beta_j$, $j=1,2$, \eqref{2.43} holds and
\eqref{2.47} now implies that
$\sigma(\ell,k_0,z)\Im(\beta_1\beta_2^*)> 0$ for $k_0$ and $\ell$
satisfying Hypothesis \ref{h2.6}, and for $z\in\bbC\backslash\bbR$.
For this choice, $\beta$ does not satisfy \eqref{BDd}. However,  we
may  normalize the boundary data as described in the proof of
Lemma \ref{l2.11}.
\end{proof}

Note that if $M\in\calD_\pm(z,k_0,\ta)$, then as a result of
Lemma \ref{l2.12} and \eqref{2.321} there is a $\beta\in\bbC^{m\times 2m}$
satisfying \eqref{2.27a} such that
\begin{equation}
M_\pm(z,k_0,\ta)=M\big(z,\ell,k_0,\ta,\tb\big).
\end{equation}
\begin{definition}
When $\calD_+(z,k_0,\ta)$ (resp., $\calD_-(z,k_0,\ta)$) is a singleton
matrix, the system \eqref{HSa} is said to be in the {\it limit point}
(l.p.) case at $\infty$ (resp., $-\infty$). If $\calD_+(z,k_0,\ta)$
(resp., $\calD_-(z,k_0,\ta)$) has nonempty interior, then \eqref{HSa} is
said to be in the {\it limit circle} (l.c.) case at $\infty$ (resp.,
$-\infty$).
\end{definition}
Indeed, for the case $m=1$, the limit point case
corresponds to a point in $\bbC$, whereas the limit circle case
corresponds to $\calD_\pm(z,k_0,\ta)$ being a closed disk in $\bbC$.

By analogy with the continuous case, these apparent geometric properties
for the disk
correspond to  analytic properties for the solutions of the
Hamiltonian system \eqref{HSa}.  To describe this correspondence, we
introduce the following spaces in which we assume that $I\subseteq\bbZ$:
\begin{subequations}\lb{2.29}
\begin{align}
\ell_A^2(I)&=\bigg\{\phi:I\to\bbC^{2m} \, \bigg| \, \sum_{k\in I}
(\phi(k),A\phi(k))_{\bbC^{2m}}<\infty \bigg\}, \lb{2.29a}
\\
N(z,\infty)&=\{\phi\in \ell_A^2([k_0,\infty)) \, | \, \calS_\rho\phi
=(zA+B)\phi \}, \lb{2.29b}
\\
N(z,-\infty)&=\{\phi\in \ell_A^2((-\infty,k_0]) \, | \,
\calS_\rho\phi=(zA+B)\phi \}, \lb{2.29c}
\end{align}
\end{subequations}
for $z\in\bbC$ and some $k_0\in\bbZ$. (Here
$(\phi,\psi)_{\bbC^n}=\sum_{j=1}^n \overline\phi_j\psi_j$
denotes the standard scalar product in $\bbC^n$, abbreviating
$\chi\in\bbC^n$ by
$\chi=(\chi_1,\dots,\chi_n)^t$.)  Both dimensions of the
spaces in \eqref{2.29b} and \eqref{2.29c},
$\dim_\bbC(N(z,\infty))$ and $\dim_\bbC(N(z,-\infty))$, are constant for
$z\in\bbC_\pm$ (see, e.g., \cite{At64}, \cite{KR74}), where
\begin{equation}
\bbC_\pm=\{\zeta\in\bbC \, | \, \pm\Im(\zeta)> 0 \}. \lb{Cpm}
\end{equation}
One then observes that  the Hamiltonian
system \eqref{HSa} is in the limit point case at
$\infty$ (resp., $-\infty$) whenever
\begin{equation}\lb{2.30}
\dim_\bbC(N(z,\infty))=m \text{ (resp., $\dim_\bbC(N(z,-\infty))=m$ 
for all $z\in\bbC\backslash\bbR$} 
\end{equation}
and in the limit circle case at $\infty$ (resp., $-\infty$) whenever
\begin{equation}\lb{2.31}
\dim_\bbC(N(z,\infty))=2m \text{ (resp., $\dim_\bbC(N(z,-\infty))=2m$)
for all $z\in\bbC$.}
\end{equation}

For the boundary condition given by
\begin{equation}
\ta\hatt{\varPsi}(z,k_0)=0
\end{equation}
with $\alpha\in\bbC^{m\times 2m}$ satisfying \eqref{r2.4a}, there
is an associated boundary set $\partial\calD_\pm(z,k_0,\ta)$ for the limiting
disk.  In either the limit point or limit circle cases, $M_\pm(z,k_0,\ta)\in
\partial\calD_\pm(a,k_0,\ta)$ is said to be a \emph{half-line
Weyl--Titchmarsh matrix}. Each such matrix is associated with the
construction of a Green's matrix for  certain boundary value problems
involving separated boundary conditions which are posed on the 
whole-line and on half-lines as will be discussed in Section \ref{s3}.

Given the definition of $M_\pm (z,k_0,\ta )$ as the limit of a sequence
$M\big(z,\ell_n,k_0,\ta,\tb\big)$, $n\in\bbN$, and given the geometry of
Weyl disks, as for the case of Hamiltonian systems of differential
equations discussed in \cite{CG01} and \cite{CG02}, we see that
$M_\pm(z,k_0,\ta)$ possesses certain properties as a complex,
matrix-valued function of $z$. For the convenience of the reader, we now
summarize some of the principal properties of half-line Weyl--Titchmarsh
matrices:
\begin{theorem}
[\cite{AD56}, \cite{Ca76}, \cite{GT00}, \cite{HS81},
\cite{HS82}, \cite{HS86}, \cite{KS88}] \lb{t2.3}
Assume Hypotheses \ref{h2.4} and let
$z\in\bbC\backslash\bbR$,
$k_0\in\bbR$, and denote by $\alpha, \gamma\in\bbC^{m\times 2m}$ matrices
satisfying \eqref{r2.4a}. Then, \\
$(i)$ $\pm M_{\pm}(z,k_0,\alpha)$ is an $m\times m$
  matrix-valued Herglotz
function of maximal rank.
In particular,
\begin{gather}
\Im(\pm M_{\pm}(z,k_0,\alpha)) > 0, \quad z\in\bbC_+, \\
M_{\pm}(\overline z,k_0,\alpha)=M_{\pm}(z,k_0,\alpha)^*, \lb{2.38} \\
\rank (M_{\pm}(z,k_0,\alpha))=m,  \\
\lim_{\varepsilon\downarrow 0} M_{\pm}(\lambda+
i\varepsilon,k_0,\alpha) \text{
exists for a.e.\
$\lambda\in\bbR$},\\
\begin{split}\lb{2.41}
M_\pm(z,k_0,\alpha) &= [-\alpha J \gamma^* +
\alpha\gamma^* M_\pm(z,k_0,\gamma)]\times \\
&\quad \times[ \alpha\gamma^*
+  \alpha J \gamma^*M_\pm(z,k_0,\gamma)]^{-1}.
\end{split}
\end{gather}
Local singularities of $\pm M_{\pm}(z,k_0,\alpha)$ and
$\mp M_{\pm}(z,k_0,\alpha)^{-1}$ are necessarily real and at most of first
order in the sense that
\begin{align}
&\mp \lim_{\epsilon\downarrow0}
\left(i\epsilon\,
M_{\pm}(\lambda+i\epsilon,k_0,\alpha)\right) \geq 0, \quad \lambda\in\bbR,
\lb{2.24b} \\
& \pm \lim_{\epsilon\downarrow0}
\left(\f{i\epsilon}{M_{\pm}(\lambda+i\epsilon,k_0,\alpha)}\right)
\geq 0, \quad \lambda\in\bbR. \lb{2.24c}
\end{align}
$(ii)$  $\pm M_{\pm}(z,k_0,\alpha)$ admit the representations
\begin{align}
&\pm M_{\pm}(z,k_0,\alpha)=F_\pm(k_0,\alpha)+\int_\bbR
d\Omega_\pm(\lambda,k_0,\alpha) \,
\big((\lambda-z)^{-1}-\lambda(1+\lambda^2)^{-1}\big) \lb{2.42} \\
&\quad =\exp\bigg(C_\pm(k_0,\alpha)+\int_\bbR d\lambda \, \Xi_{\pm}
(\lambda,k_0,\alpha)
\big((\lambda-z)^{-1}-\lambda(1+\lambda^2)^{-1}\big)
\bigg), \lb{2.430}
\end{align}
where
\begin{align}
F_\pm(k_0,\alpha)&=F_\pm(k_0,\alpha)^*, \quad \int_\bbR
\|d\Omega_\pm(\lambda,k_0,\alpha)\|_{\bbC^{m\times m}} \,
(1+\lambda^2)^{-1}<\infty,
  \\
C_\pm(k_0,\alpha)&=C_\pm(k_0,\alpha)^*,
\quad 0\le\Xi_\pm(\dott,k_0,\alpha)
\le I_m \, \text{  a.e.}
\end{align}
Moreover,
\begin{align}
\Omega_\pm((\lambda,\mu],k_0,\alpha)&
=\lim_{\delta\downarrow
0}\lim_{\varepsilon\downarrow 0}\f1\pi
\int_{\lambda+\delta}^{\mu+\delta} d\nu \, \Im(\pm
M_\pm(\nu+i\varepsilon,k_0,\alpha)),  \\
\Xi_\pm(\lambda,k_0,\alpha)&=\lim_{\varepsilon\downarrow 0}
\pi^{-1}\Im(\ln(\pm
M_\pm(\lambda+i\varepsilon,k_0,\alpha))) \text{ for a.e.\
$\lambda\in\bbR$}.
\end{align}
$(iii)$  Define the $2m\times m$ matrices
\begin{align}
U_\pm(z,\ell,k_0,\alpha)&=\begin{pmatrix}u_{\pm,1}(z,\ell,k_0,\alpha) \\
u_{\pm,2}(z,\ell,k_0,\alpha)  \end{pmatrix}
=\Psi(z,\ell,k_0,\alpha)\begin{pmatrix} I_m \\
M_\pm(z,k_0,\alpha) \end{pmatrix}  \no \\
&=\begin{pmatrix}\theta_1(z,\ell,k_0,\alpha)
& \phi_1(z,\ell,k_0,\alpha)\\
\theta_2(z,\ell,k_0,\alpha)
& \phi_2(z,\ell,k_0,\alpha)\end{pmatrix}
\begin{pmatrix} I_m \\
M_\pm(z,k_0,\alpha) \end{pmatrix}, \lb{2.52}
\end{align}
with $\theta_j(z,\ell,k_0,\alpha)$, and
$\phi_j(z,\ell,k_0,\alpha)$, $j=1,2$, defined by \eqref{FSc}. Then, for
every $\zeta\in\bbC^{2m}$, 
$U_\pm(z,\ell,k_0,\alpha)\zeta\in\ell_A^2([k,\pm\infty))$.
Moreover,
\begin{align}
\Im(M_+(z,k_0,\alpha))&=\Im(z) \sum_{^+[k_0,\infty)}
U_+(z,s,k_0,\alpha)^* A(s)
U_+(z,s,k_0,\alpha), \label{2.88}\\
\Im(M_-(z,k_0,\alpha))&=\Im(z) \sum_{(-\infty,k_0]}
U_-(z,s,k_0,\alpha)^* A(s)
U_-(z,s,k_0,\alpha).\label{2.89}
\end{align}
\end{theorem}

We conclude this section by first giving another characterization of
the elements of the limiting disks $\calD_\pm(z,k_0,\ta)$ and then noting
a connection between these elements and certain solutions
of a related Riccati equation.

\begin{lemma}\label{l3.3}
Assume Hypothesis \ref{h2.6} and let $z\in\bbC\backslash\bbR$. Moreover, 
suppose that $U(z,k,k_0,\ta)$ is defined by \eqref{2.19} in terms of an
$M\in\bbC^{m\times m}$ so that the columns of $U(z,k,k_0,\ta)$ are in
$\ell_A^2([k_0,\pm\infty))$. Then, $M\in\calD_\pm(z,k_0,\ta)$  if and
only if for $k\in\bbZ$,
\begin{equation}\label{3.3}
-\sigma(k,k_0,z)\Im (u_1(z,k)^*\rho(k)u_2^+(z,k)) > 0,
\end{equation}
or equivalently,
\begin{equation}\label{3.4}
-\sigma(k,k_0,z)\Im (\rho(k)u_2^+(z,k)u_1(z,k)^{-1}) >0.
\end{equation}
\end{lemma}
\begin{proof}
By \eqref{2.21} and \eqref{2.32a},
\begin{align}
-2\sigma(k,k_0,z)\Im (u_1(z,k)^*\rho(k)u_2^+(z,k)) &=
\sigma(k,k_0,z)\hatt{U}^*(iJ_{\rho})\hatt{U}{\big |}_k \notag\\
&= -E_\ell (M) + 2|\Im(z)|\sum_{^+[k,\ell]} U^*A U.\label{2.920}
\end{align}
The summation expression in \eqref{2.920} is positive by
Hypothesis \ref{h2.6}, and is decreasing as a function of $k$ with fixed
$\ell$ while increasing as a function of $\ell$ with fixed $k$. As a
consequence,  should a column of $U(z,k,k_0,\ta)$ not be in
$\ell_A^2([k_0,\pm\infty))$, then  there is a vector
$\zeta\in\bbC^{2m}$ such that $\zeta^*E_\ell(M)\zeta>0$ for 
large $|\ell|$ and hence $M\not\in\calD_\pm(z,k_0,\ta)$. Thus, we assume
that the columns of $U(z,k,k_0,\ta)$ are in $\ell_A^2([k_0,\pm\infty))$.

If $M\in\calD(z,k_0,\ta)$ then \eqref{3.3} follows because 
$-E_\ell(M)\ge 0$
and because the sum in \eqref{2.920} is positive for large $|\ell|$. If
$M\not\in\calD(z,k_0,\ta)$ then $\zeta^*E_{\ell_0}(M)\zeta>0$ for some
$\zeta\in\bbC^{2m}$ and some $\ell_0$ and hence
$\zeta^*E_{\ell}(M)\zeta\ge \zeta^*E_{\ell_0}(M)\zeta$
for $|\ell| > |\ell_0|$. Then, for sufficiently large $|\ell|$ and 
$|k|$, the left-hand side of \eqref{3.3} is negative.

The equivalence of \eqref{3.3} and \eqref{3.4} follows because
\begin{equation}\label{3.5}
u_1^*\Im (\rho u_2^+u_1^{-1})u_1=\Im (u_1^*\rho u_2^+)>0.
\end{equation}
\end{proof}

The Hamiltonian system \eqref{HS}, described in Definition \ref{d2.4},
can be written as
\begin{subequations}
\begin{align}
\rho\psi^+_2 &= (zA_{1,1}+B_{1,1})\psi_1 + (zA_{1,2}+B_{1,2})\psi_2,
\label{2.73a}\\
\rho^-\psi^-_1 &= (zA_{2,1}+B_{2,1})\psi_1 +
(zA_{2,2}+B_{2,2})\psi_2.\label{2.73b}
\end{align}
\end{subequations}
When $\varPsi(z,k)\in\mathbb{C}^{2m\times m}$  represents a solution
of \eqref{HSa}
for which $\psi_j(z,k)\in\mathbb{C}^{m\times m}$, $j=1,2$,
are nonsingular for $k\in [k_0,\ell]$, then \eqref{2.73a} and \eqref{2.73b}
respectively yield,
\begin{subequations}
\begin{align}
\rho\psi^+_2\psi^{-1}_1 &= (zA_{1,1}+B_{1,1}) +
(zA_{1,2}+B_{1,2})\psi_2\psi^{-1}_1,\\
\psi^{-1}_1&= [\rho^-\psi^-_1 -
(zA_{2,2}+B_{2,2})\psi_2]^{-1}(zA_{2,1}+B_{2,1}),
\end{align}
\end{subequations}
from which it follows that
\begin{equation}
\begin{split}
\rho\psi^+_2\psi^{-1}_1&=(zA_{1,1}+B_{1,1}) + (zA_{1,2}+B_{1,2})
[\rho^-\psi^-_1\psi_2^{-1} - (zA_{2,2}+B_{2,2})]^{-1}\\
&\quad\times(zA_{2,1}+B_{2,1}).
\end{split}
\end{equation}
Thus, $V(z,k)=\rho(k)\psi^+_2(z,k)\psi_1(z,k)^{-1}$, $k\in[k_0,\ell]$,
yields a solution
of the Riccati equation given by
\begin{equation}\label{Ric}
\begin{split}
V&=(zA_{1,1}+B_{1,1}) + (zA_{1,2}+B_{1,2})
[\rho^-(V^-)^{-1}\rho^- - (zA_{2,2}+B_{2,2})]^{-1}\\
&\quad\times(zA_{2,1}+B_{2,1}).
\end{split}
\end{equation}
As an immediate consequence of these observations, we obtain the 
following result.

\begin{lemma}\label{l2.19}
Assume Hypothesis \ref{h2.6} and let $z\in\bbC\backslash\bbR$. Moreover,
suppose that $U(z,k,k_0,\ta)$ is defined by \eqref{2.19} in terms of
an $M\in\bbC^{m\times m}$ so that the columns of $U(z,k,k_0,\ta)$ are in
$\ell_A^2([k_0,\pm\infty))$. Then, $M\in\calD_\pm(z,k_0,\ta)$  if and
only if for $k\in\bbZ$,
$V(z,k)=\rho(k)u_2^+(z,k)u_1(z,k)^{-1}$ represents
a solution of the Riccati equation \eqref{Ric} for $z\in\bbC\backslash\bbR$.
\end{lemma}

\section{Boundary Value Problems and Green Matrices on the Whole Line and
on Half-Lines}
\label{s3}

In this section, we consider the nonhomogeneous equation given by
\begin{equation}\label{3.1}
\calS_\rho\psi=(zA+B)\psi +Af
\end{equation}
associated with the Hamiltonian system \eqref{HSa} on the whole-line or
on half-lines.

First, we discuss the whole-line case and assume Hypothesis
\ref{h2.4} on $\bbZ$. Hence, we consider \eqref{3.1} with
$f\in\ell_A^2(\bbZ)$. To keep matters reasonably short, the endpoints
$\infty$ and
$-\infty$ will separately be assumed to be either of \emph{limit point}
or \emph{limit circle} type for
\eqref{HSa}. (These cases typically receive most attention and the sequel
\cite{CGR04} to this paper will, in particular, focus on the limit point
case at $\infty$ and $-\infty$.) We describe for $k,\ell\in\bbZ$, and
$z\in\bbC\backslash\bbR$, a matrix
$K(z,k,\ell)\in\bbC^{2m\times 2m}$ for which the following properties hold:
\begin{equation}\label{3.2}
\sum_{\ell\in \bbZ}K(z,k,\ell)A(\ell)K(z,k,\ell)^*<\infty, \quad k\in\bbZ.
\end{equation}
If $f\in\ell_A^2(\bbZ)$ and if
\begin{equation}\label{3.30}
y(z,k)=\sum_{\ell\in \bbZ}K(z,k,\ell)A(\ell)f(\ell), \quad k\in\bbZ,
\end{equation}
then $y(z,\cdot)\in\ell_A^2(\bbZ)$ and $y(z,\cdot)$ satisfies \eqref{3.1} on
$\bbZ$.  In addition, it will be seen that $y(z,\cdot)$ satisfies
certain boundary conditions at $\infty$ (resp., $-\infty$) if \eqref{HSa}
is in the limit circle case at $\infty$ (resp., $-\infty$).

As a matter of convenience, we state the next
theorem assuming that $z\in\bbC_+$ (cf.\ \eqref{Cpm}) and note that the
theorem can be restated for $z\in\bbC_-$ with the details of the proof
essentially unchanged.
\begin{theorem}\label{t2.20}
Assume Hypothesis \ref{h2.4} on $\bbZ$ and suppose that $z\in\bbC_+$ and 
$k,\ell\in\bbZ$. Let
$K(z,k,\ell)$ be defined by
\begin{equation}\label{2.80}
K(z,k,\ell)=
\begin{cases}
U_+(z,k)\omega(z) U_-(\bar z,\ell)^*, & k>\ell, \\
\begin{pmatrix}
u_{+,1}(z,k)\omega(z) u_{-,1}(\bar z,k)^*&
u_{+,1}(z,k)\omega(z) u_{-,2}(\bar z,k)^*\\
u_{-,2}(z,k)\omega(z) u_{+,1}(\bar z,k)^*&
u_{-,2}(z,k)\omega(z) u_{+,2}(\bar z,k)^*
\end{pmatrix}, & k=\ell, \\[1mm]
U_-(z,k)\omega(z) U_+(\bar z,\ell)^*, & k<\ell.  
\end{cases}
\end{equation} 
Here $U_{\pm}(z,k)$ is defined in \eqref{2.19} with
$M=M_{\pm}(z)\in\partial\calD_{\pm}(z,k_0,\ta)$,
$U_{\pm}(\bar z,k)$ is defined in \eqref{2.19} with
$M=M_\pm(\bar z)=M_{\pm}(z)^*\in\partial\calD_{\pm}(\bar z,k_0,\ta)$, and
\begin{equation}
U_{\pm}(z,k)= \begin{pmatrix}u_{\pm,1}(z,k)\\ u_{\pm,2}(z,k)
\end{pmatrix}
\end{equation}
with $u_{\pm,j}(z,k)\in\bbC^{m\times m}$, $j=1,2$. Moreover,
\begin{equation}\label{2.82}
\omega(z)=\pm\big(\hatt{U}_{\mp}(\bar z,
k_0)^*J_{\rho}(k_0)\hatt{U}_{\pm}(z,k_0)\big)^{-1}=[M_-(z) - M_+(z)]^{-1}.
\end{equation}
With $K(z,k,\ell)$ so defined, \eqref{3.2} is satisfied. Moreover, as 
defined in \eqref{3.30},
$y(z,\cdot)$ satisfies \eqref{3.1} on $\bbZ$ and is in $\ell_A^2(\bbZ)$.
\end{theorem}
\begin{proof}
We begin by defining notation to be used for the remainder of this
section. We adopt the following convention:
\begin{equation}\lb{2.84}
F^{\circledast}(z,k)=F(\bar z,k)^* \, \text{ for $F\in\bbC^{m\times
m}$ (or $F\in\bbC^{2m\times m}$)}.
\end{equation}
Let the matrices $a(z,k)$, $b(z,k)$, $c(z,k)\in\bbC^{m\times m}$ be defined by
\begin{equation}
zA(k) + B(k)= \begin{pmatrix}a(z,k)& b(z,k)\\ b(\bar z, k)^* & c(z,k)
\end{pmatrix}=
\begin{pmatrix}a(z,k)& b(z,k)\\ b^{\circledast}( z, k) & c(z,k) \end{pmatrix},
\end{equation}
where $A(k)$, $B(k)$ are given in \eqref{d2.1b} and \eqref{d2.1c}
respectively, and are subject
to Hypothesis \ref{h2.4} on $\bbZ$. We also note that
$a^{\circledast}(z,k)=a(z,k)$, and that
$c^{\circledast}(z,k)=c(z,k)$. Lastly, let $\vp_{\pm}(z,k)$ and
$\vt_{\pm}(z,k)$ be defined by
\begin{equation}\label{2.86}
\vp_{\pm}(z,k)=u_{\pm,1}(z,k),\qquad \vt_{\pm}(z,k)=u_{\pm,2}(z,k).
\end{equation}

With this convention,
$\hatt{U}_{\mp}^{\circledast}(z, k)J_{\rho}(k)\hatt{U}_{\pm}(z,k)=
\hatt{U}_{\mp}(\bar z, k)^*J_{\rho}(k)\hatt{U}_{\pm}(z,k)$. Then,
  by \eqref{2.21} we note that
\begin{subequations}\label{2.83}
\begin{align}\label{2.83a}
\hatt{U}_{\mp}^{\circledast}(z, k)J_{\rho}(k)\hatt{U}_{\pm}(z,k)&=
\hatt{U}_{\mp}^{\circledast}(z, k_0)J_{\rho}(k_0)\hatt{U}_{\pm}(z,k_0)\\
&=M_{\mp}(z) - M_{\pm}(z).
\end{align}
\end{subequations}
Given that $\Im(M_{\pm}(z))\gtrless 0$, we see that
$\hatt{U}_{\mp}^{\circledast}(z, k)J_{\rho}(k)\hatt{U}_{\pm}(z,k)$ is
invertible for $k\in\bbZ$.

Next, we note that for $k\ne\ell$, $0=\big((\calS_{\rho} - zA -
B)K(z,\cdot,\ell)\big)(k)$.
Thus, to verify that $y$ defined in \eqref{3.30} solves \eqref{3.1}, it is necessary
to show that
\begin{subequations}\label{2.87}
\begin{align}
I_{2m}&=\big((\calS_{\rho} - zA - B)K(z,\cdot,k)\big)(k)\\
&=\big(\calS_{\rho}K(z,\cdot,k)\big)(k) - \begin{pmatrix}a(z,k)& b(z,k)\\
b^{\circledast}( z, k) & c(z,k) \end{pmatrix}K(z,k,k).
\end{align}
\end{subequations}
Then, by \eqref{2.80} for $k=\ell$, and \eqref{2.86},
\begin{equation}
\begin{pmatrix}a(z,k)& b(z,k)\\
b^{\circledast}( z, k) & c(z,k) \end{pmatrix}K(z,k,k)=
\begin{pmatrix}a& b\\
b^{\circledast} & c \end{pmatrix}
\begin{pmatrix}\vp_+\omega\vp_-^{\circledast} &
\vp_+\omega\vt_-^{\circledast}\\[1mm]
\vt_-\omega\vp_+^{\circledast} & \vt_-\omega\vt_+^{\circledast}
\end{pmatrix},
\end{equation}
and by \eqref{1.10} and \eqref{2.80} for $k\neq\ell$,
\begin{equation}
\big(\calS_{\rho}K(z,\cdot,k)\big)(k)=\begin{pmatrix}
\rho\vt_+^+\omega\vp_-^{\circledast} &
\rho\vt_+^+\omega\vt_-^{\circledast}\\[1mm]
\rho^-\vp_-^-\omega\vp_+^{\circledast} &
\rho^-\vp_-^-\omega\vt_+^{\circledast}
\end{pmatrix}.
\end{equation}
Thus, \eqref{2.87} is equivalent to the following system:
\begin{subequations}\label{2.90}
\begin{align}
I_m &= \rho\vt_+^+\omega\vp_-^{\circledast}
-a\vp_+\omega\vp_-^{\circledast}
-b\vt_-\omega\vp_+^{\circledast},\\
0_m &= \rho\vt_+^+\omega\vt_-^{\circledast}
-a\vp_+\omega\vt_-^{\circledast}
-b\vt_-\omega\vt_+^{\circledast},\\
0_m &= \rho^-\vp_-^-\omega\vp_+^{\circledast}
-b^{\circledast}\vp_+\omega\vp_-^{\circledast} -
c\vt_-\omega\vp_+^{\circledast},\\
I_m &= \rho^-\vp_-^-\omega\vt_+^{\circledast}
-b^{\circledast}\vp_+\omega\vt_-^{\circledast} -
c\vt_-\omega\vt_+^{\circledast}.
\end{align}
\end{subequations}
Given that  $0=\big((\calS_{\rho} - zA - B)U_{\pm}(z,\cdot)\big)(k)$,
we obtain additionally that
\begin{subequations}\label{2.91}
\begin{align}
\rho\vt_{\pm}^+   &= a\vp_{\pm} + b\vt_{\pm},\label{2.91a}\\
\rho^-\vp_{\pm}^- &= b^{\circledast}\vp_{\pm} + c\vt_{\pm}.\label{2.91b}
\end{align}
\end{subequations}
 From \eqref{2.91}, we obtain
\begin{subequations}\label{2.92}
\begin{align}
(\vt_{\pm}^{\circledast})^{+}\rho &= \vp_{\pm}^{\circledast}a
+ \vt_{\pm}^{\circledast}b^{\circledast},\label{2.92a}\\
(\vp_{\pm}^{\circledast})^{-}\rho^- &= \vp_{\pm}^{\circledast}b
+ \vt_{\pm}^{\circledast}c.\label{2.92b}
\end{align}
\end{subequations}

Substituting into \eqref{2.90} the expressions for $\rho\vt_+^+$ and for
$\rho_-\vp_-^-$ in \eqref{2.91},
we obtain the equivalent system
\begin{subequations}\label{2.93}
\begin{align}
I_m &= b\vt_+\omega\vp_-^{\circledast}
-b\vt_-\omega\vp_+^{\circledast},\label{2.93a}\\
0_m &= b\vt_+\omega\vt_-^{\circledast}
-b\vt_-\omega\vt_+^{\circledast},\label{2.93b}\\
0_m &= b^{\circledast}\vp_-\omega\vp_+^{\circledast} -
b^{\circledast}\vp_+\omega\vp_-^{\circledast},\label{2.93c}\\
I_m &= b^{\circledast}\vp_-\omega\vt_+^{\circledast} -
b^{\circledast}\vp_+\omega\vt_-^{\circledast}.\label{2.93d}
\end{align}
\end{subequations}
Verification of \eqref{2.87} comes with first showing the consistency of
the equations in \eqref{2.93}
and hence the consistency of those in \eqref{2.90}. Following this, we identify
$\omega$ by \eqref{2.82} as terms by which the equations in
\eqref{2.93} hold.

Given that $b(z,k)$ and $b(z,k)^\circledast$ are invertible by
Hypothesis \ref{h2.1}, and given the invertibility of $\vp_{\pm}(z,k)$
and $\vt_{\pm}(z,k)$, $z\in\bbC_+$, $k\in\bbZ$ which follows from
Lemma \ref{l2.19}, we obtain from \eqref{2.93b}
\begin{equation}\label{2.94}
\omega = \vt_+^{-1}\vt_-\omega\vt_+^{\circledast}
(\vt_-^{\circledast})^{-1},
\end{equation}
and from \eqref{2.93c}
\begin{equation}\label{2.95}
\omega = \vp_+^{-1}\vp_-\omega\vp_+^{\circledast}
(\vp_-^{\circledast})^{-1}.
\end{equation}
Replacing in \eqref{2.93a} the expression for $\omega$ in
\eqref{2.95} one obtains
\begin{equation}\label{2.96}
\omega = \vp_-^{-1}(\vt_+\vp_+^{-1} - \vt_-\vp_-^{-1})^{-1}b^{-1}
(\vp_+^{\circledast})^{-1},
\end{equation}
and replacing in \eqref{2.93d} the expression for $\omega$ in \eqref{2.95}
one obtains
\begin{equation}\label{2.97}
\omega = \vp_-^{-1}(b^{\circledast})^{-1}\big((\vp_+^{\circledast})^{-1}
\vt_+^{\circledast} - (\vp_-^{\circledast})^{-1}
\vt_-^{\circledast}\big)^{-1}(\vp_+^{\circledast})^{-1}.
\end{equation}
Equating the right-hand sides of \eqref{2.96} and \eqref{2.97} then yields
\begin{equation}\label{2.98}
b(\vt_+\vp_+^{-1} - \vt_-\vp_-^{-1}) = \big((\vp_+^{\circledast})^{-1}
\vt_+^{\circledast} - (\vp_-^{\circledast})^{-1}
\vt_-^{\circledast}\big)b^{\circledast}.
\end{equation}
Substituting on the left-hand side of \eqref{2.98} using the expression
for $b\vt_{\pm}$ given in \eqref{2.91a}, and
substituting on the right-hand side of \eqref{2.98} using the expression
for $\vt_{\pm}^{\circledast}b^{\circledast}$ given
in \eqref{2.92a} one obtains
\begin{equation}\label{2.100}
\rho\vt_+^+\vp_+^{-1} - \rho\vt_-^+\vp_-^{-1} =
(\vp_+^{\circledast})^{-1}(\vt_+^{\circledast})^{+}\rho -
(\vp_-^{\circledast})^{-1}(\vt_-^{\circledast})^{+}\rho.
\end{equation}
That equation \eqref{2.100} holds follows from \eqref{2.21}, from
the fact that
\begin{equation}\label{2.101}
\hatt{U}_{\pm}^{\circledast}(z, k)J_{\rho}(k)\hatt{U}_{\pm}(z,k)=
\hatt{U}_{\pm}^{\circledast}(z,
k_0)J_{\rho}(k_0)\hatt{U}_{\pm}(z,k_0)=0, \\
\end{equation}
and hence that
\begin{equation}\label{2.102}
(\vt_{\pm}^{\circledast})^{+}\rho\vp_{\pm} =
\vp_{\pm}^{\circledast}\rho\vt_{\pm}^+.
\end{equation}
As a consequence, we see that \eqref{2.93a}, \eqref{2.93c},
\eqref{2.93d} are consistent.

Replacing $\omega$ in \eqref{2.95} with the expression for
$\omega$ in \eqref{2.96} yields
\begin{equation}\label{2.103}
\omega = \vp_+^{-1}(\vt_+\vp_+^{-1} - \vt_-\vp_-^{-1})^{-1}
b^{-1}(\vp_-^{\circledast})^{-1},
\end{equation}
and replacing $\omega$ in \eqref{2.94} with the expression for
$\omega$ in \eqref{2.96} yields
\begin{equation}\label{2.104}
\omega = \big( \vt_-^{\circledast}(\vt_+^{\circledast})^{-1}
\vp_+^{\circledast}b(\vt_+\vp_+^{-1}\vp_-\vt_-^{-1}\vt_+ -
\vt_+)\big)^{-1}.
\end{equation}
Showing the equivalence of \eqref{2.103} and \eqref{2.104} will yield
the equivalence of \eqref{2.93b} and \eqref{2.93c},
and hence the consistency of all equations in \eqref{2.93}.

Replacing $\vp_+^{\circledast}b$ in \eqref{2.104} with the expression
obtained from \eqref{2.92b}
yields
\begin{align}
\omega &= \big[ \vt_-^{\circledast}((\vt_+^{\circledast})^{-1}
(\vp_+^{\circledast})^{-}\rho^- - c)(\vt_+\vp_+^{-1}\vp_-\vt_-^{-1}\vt_+
- \vt_+)\big]^{-1}.\\
\intertext{By \eqref{2.102}, $(\vt_+^{\circledast})^{-1}
(\vp_+^{\circledast})^{-}\rho^-  = \rho^-\vp_+^-\vt_+^{-1}$, and thus,}
\omega&= \big[\vt_-^{\circledast}(\rho^-\vp_+^-\vt_+^{-1}-c)
(\vt_+\vp_+^{-1}\vp_-\vt_-^{-1}\vt_+-\vt_+) \big]^{-1}.\\
\intertext{By \eqref{2.91b}, $\rho^-\vp_+^-\vt_+^{-1}-c
= b^{\circledast}\vp_+\vt_+^{-1}$, and hence,}
\omega&= \big[ \vt_-^{\circledast}b^{\circledast}
\vp_+\vt_+^{-1}(\vt_+\vp_+^{-1}\vp_-\vt_-^{-1} - I_m)\vt_+  \big]^{-1},\\
&= \big[\vt_-^{\circledast}b^{\circledast}(\vp_-\vt_-^{-1}
- \vp_+\vt_+^{-1})\vt_+  \big]^{-1},\\
&= \vp_+^{-1}\big[ (\vp_-^{\circledast})^{-1}\vt_-^{\circledast}
b^{\circledast}(\vp_-\vt_-^{-1}- \vp_+\vt_+^{-1})
\vt_+\vp_+^{-1} \big]^{-1}(\vp_-^{\circledast})^{-1}.\\
\intertext{By \eqref{2.92a}, $\vt_-^{\circledast}b^{\circledast}
=(\vt_-^{\circledast})^{+}\rho - \vp_-^{\circledast}a$, and thus,}
\omega&= \vp_+^{-1}\big[ ((\vp_-^{\circledast})^{-1}
(\vt_-^{\circledast})^{+}\rho
- a)  (\vp_-\vt_-^{-1}\vt_+\vp_+^{-1} - I_m)\big]^{-1}
(\vp_-^{\circledast})^{-1}.\\
\intertext{By \eqref{2.102}, $(\vp_-^{\circledast})^{-1}
(\vt_-^{\circledast})^{+}\rho  = \rho\vt_-^+\vp_-^{-1}$, and hence,}
\omega&= \vp_+^{-1}\big[(\rho\vt_-^+\vp_-^{-1} - a)
(\vp_-\vt_-^{-1}\vt_+\vp_+^{-1} - I_m) 
\big]^{-1}(\vp_-^{\circledast})^{-1}\\ 
&= \vp_+^{-1}\big[(\rho\vt_-^+ -
a\vp_-)(\vt_-^{-1}\vt_+\vp_+^{-1} - \vp_-^{-1}) 
\big]^{-1}(\vp_-^{\circledast})^{-1}. \label{2.112}
\end{align}
By \eqref{2.91a}, $\rho\vt_-^+ - a\vp_- = b\vt_-$. Using this equivalence in
\eqref{2.112} yields the right-hand side of \eqref{2.103} and thus
establishes the equivalence of \eqref{2.103} and \eqref{2.104}. Thus the
equations in \eqref{2.90} and \eqref{2.93} are consistent.

To obtain \eqref{2.82}, we first note that \eqref{2.96} and \eqref{2.103}
can be written as
\begin{align}
\omega &= \vp_{\pm}^{-1}(\vt_+\vp_+^{-1}-\vt_-\vp_-^{-1})^{-1}b^{-1}
(\vp_{\mp}^{\circledast})^{-1}\\
&= \big[ \vp_{\mp}^{\circledast} (b \vt_+\vp_+^{-1}-b \vt_-\vp_-^{-1})
\vp_{\pm}  \big]^{-1}. \\
\intertext{By \eqref{2.91a}, $b\vt_{\pm}=\rho\vt_{\pm}^+ - a\vp_{\pm}$,
and hence,}
\omega &= \big[ \vp_{\mp}^{\circledast} (\rho \vt_{+}^+\vp_{+}^{-1}
- \rho \vt_{-}^+\vp_{-}^{-1})\vp_{\pm}  \big]^{-1}.\label{2.115}\\
&= \pm\big[ \vp_{\mp}^{\circledast} (\rho \vt_{\pm}^+\vp_{\pm}^{-1}
- \rho \vt_{\mp}^+\vp_{\mp}^{-1})\vp_{\pm}  \big]^{-1}.\\
\intertext{By \eqref{2.102}, $ \rho\vt_{\mp}^+\vp_{\mp}^{-1}
= (\vp_{\mp}^{\circledast})^{-1}(\vt_{\mp}^{\circledast})^{+}\rho$, and
thus,}
\omega &= \pm\big[ \vp_{\mp}^{\circledast}\rho \vt_{\pm}^+
- (\vt_{\mp}^{\circledast})^{+}\rho \vp_{\pm}  \big]^{-1}
= \pm \big(\hatt{U}_{\mp}^{\circledast}J_{\rho}\hatt{U}_{\pm}\big)^{-1}
\end{align}
which by \eqref{2.83a} yields \eqref{2.82}, and thus completes the demonstration that
$y$, as defined by \eqref{3.30}, satisfies \eqref{3.1}.

For the remaining properties of $K(z,k,\ell)$, we first note that \eqref{2.88} 
and \eqref{2.89} 
imply that $K(z,k,\ell)$ satisfies \eqref{3.2}. And finally, by an argument in direct analogy
with that given  in \cite[Lemma 4.2 ]{HS83} (see also
\cite[Lemma 4.1]{HS84}), we see that $y$ given by \eqref{3.30} satisfies
\begin{equation}
\sum_{k\in\bbZ}y(z,k)^*A(k)y(z,k) \le (\Im
(z))^{-2}\sum_{k\in\bbZ}f(k)^*A(k)f(k).
\end{equation}
As a result, we note that $y(z,\cdot)\in\ell^2_A(\bbZ)$ whenever
$f\in\ell^2_A(\bbZ)$.
\end{proof}

For $k\ne\ell$, there exists an alternative expression for $K(z,k,\ell)$
in terms of the half-line M-matrices and the fundamental solution
$\Psi(z,k,k_0,\ta)$ which is defined in \eqref{FS}. In direct analogy
with equation (4.5) of \cite{HS83}, we note that 
\begin{align}
&K(z,k,\ell)=U_+(z,k)\omega(z) U_-(\bar z,\ell)^*=U_+(z,k)[M_-(z) -
M_+(z)]^{-1} U_-(\bar z,\ell)^* \no \\ 
&\quad =\Psi(z,k,k_0,\ta)\notag\\
&\qquad\times\begin{pmatrix}
[M_-(z) - M_+(z)]^{-1} & [M_-(z) - M_+(z)]^{-1}M_-(z)\\
M_+(z)[M_-(z) - M_+(z)]^{-1} & M_+(z)[M_-(z) - M_+(z)]^{-1}M_-(z)
\end{pmatrix} \no \\
&\qquad\times\Psi(\bar z,\ell ,k_0,\ta)^*, \quad k>\ell
\end{align}
and  
\begin{align}
&K(z,k,\ell)=U_-(z,k)\omega(z) U_+(\bar z,\ell)^*=U_-(z,k)[M_-(z) -
M_+(z)]^{-1} U_+(\bar z,\ell)^* \no \\ 
&\quad =\Psi(z,k,k_0,\ta)\notag\\
&\qquad\times\begin{pmatrix}
[M_-(z) - M_+(z)]^{-1} & [M_-(z) - M_+(z)]^{-1}M_+(z)\\
M_-(z)[M_-(z) - M_+(z)]^{-1} & M_-(z)[M_-(z) - M_+(z)]^{-1}M_+(z)
\end{pmatrix} \no \\
& \qquad \times\Psi(\bar z,\ell ,k_0,\ta)^*, \quad k<\ell,
\end{align}
noting that $M_+(M_- - M_+)^{-1}M_- = M_-(M_- - M_+)^{-1}M_+$.

Given the notation introduced in Theorem \ref{t2.20}, specifically in
\eqref{2.86}, let
\begin{equation}
V_{\pm}(z,k)= \rho(k)u_{\pm 2}^+(z,k)u_{\pm 1}(z,k)^{-1}
= \rho(k)\vt_{\pm}^+(z,k)\vp_{\pm}(z,k)^{-1}.
\end{equation}
We note that $V_{\pm}(z,k)$ is a solution of the Riccati equation
given in \eqref{Ric}. We also note that by \eqref{2.115}, 
\begin{align}
\omega= (M_- - M_+)^{-1}&=\vp_{\pm}^{-1}(\rho\vt_+^+\vp_+^{-1}
- \rho\vt_-^+\vp_-^{-1}   )^{-1}(\vp_{\mp}^{\circledast})^{-1}\\
&= \vp_{\pm}^{-1}(V_+ - V_-)^{-1}(\vp_{\mp}^{\circledast})^{-1}.
\end{align}
Then, by \eqref{2.80} for $k=\ell$,
\begin{equation}
K(z,k,k)=
\begin{pmatrix}\vp_+(z,k)\omega(z)\vp_-^{\circledast}(z,k)
& \vp_+(z.k)\omega(z)\vt_-^{\circledast}(z,k) \\[1mm]
\vt_-(z,k)\omega(z)\vp_+^{\circledast}(z,k)
& \vt_-(z,k)\omega(z)\vt_+^{\circledast}(z,k)
\end{pmatrix}.
\end{equation}
An alternative representation for the entries of
the matrix $K(z,k,k)$ also exists:
\begin{align}
\vp_+\omega\vp_-^{\circledast} &= (V_+ - V_-)^{-1},\\
\vp_+\omega\vt_-^{\circledast}&= (V_+ - V_-)^{-1}(\vp_-^{\circledast})^{-1}\vt_-^{\circledast},\notag\\
&=(V_+ - V_-)^{-1}(u_{-, 1}^{\circledast})^{-1}u_{-, 2}^{\circledast}, \\
\vt_-\omega\vp_+^{\circledast} &= \vt_-\vp_-^{-1}(V_+ - V_-)^{-1},\notag\\
&= u_{-, 2}(u_{-, 1})^{-1}(V_+ - V_-)^{-1}  \\
\vt_-\omega\vt_+^{\circledast} &= \vt_-\vp_-^{-1}(V_+ - V_-)^{-1}(\vp_+^{\circledast})^{-1}\vt_+^{\circledast}\notag\\
&= u_{-,2}(u_{-, 1})^{-1}(V_+ - V_-)^{-1}(u_{+,1}^{\circledast})^{-1}u_{+, 2}^{\circledast}.
\end{align}

The proof of the next lemma relies upon an argument involving the geometry 
of the Weyl disks described in Definition \ref{dWD}. We refer the reader
to \cite[Theorem 2.1]{HS84} for details while noting that the discussion
in \cite{HS84} occurs in the context of Hamiltonian systems of ordinary
differential equations but that the argument remains the same for the
current setting.
\begin{lemma}\label{l3.2}
Assume that the Hamiltonian system \eqref{HSa} which satisfies
Hypothesis \ref{h2.4} is in the limit point or the limit circle case at
$\infty$. Let $z_1,z_2\in\bbC\backslash\bbR$. Then,
\begin{equation}
\lim_{k\to\infty}\begin{bmatrix} I & M_+(z_2)^* \end{bmatrix}
\hatt{\Psi}(z,k,k_0,\ta)^*J_\rho\hatt{\Psi}(z,k,k_0,\ta)
\begin{bmatrix} I \\ M_+(z_1) \end{bmatrix}=0,
\end{equation}
where $M_+(z_j)\in\partial\calD(z_j,k_0,\ta)$, $j=1,2$, and where
$\Psi(z,k,k_0,\ta)$ is the fundamental matrix defined in \eqref{FS}.
\end{lemma}
Of course, an analogous result can be stated when the Hamiltonian system \eqref{HSa} is in the limit point
or the limit circle case at $-\infty$. Moreover, an immediate consequence of this result, like that of its
continuous counterpart, is the following corollary. Again see \cite[Corollary~2.3]{HS84} for details of the
proof that also remains the same for the current setting.
\begin{corollary}\label{c3.3}
The Hamiltonian system \eqref{HSa} which satisfies Hypothesis \ref{h2.4}
is in the limit point case at
$\infty$ if and only if
\begin{equation}
\lim_{k\to\infty} \hatt{y}(z_1,k)^*J_\rho \hatt{y}(z_2,k) = 0, \quad 
z_1,z_2\in\bbC\backslash\bbR
\end{equation}
for all $\ell_A^2([k_0,\infty))$-solutions $y(z_j,\cdot)$, $j=1,2$, of
\eqref{HSa}.
\end{corollary}
As before, we note that an analogous result for Corollary~\ref{c3.3} can be stated when the
Hamiltonian system \eqref{HSa} is in the limit point case at $-\infty$.

As a consequence of the preceding results of this section, we have the 
following theorem which effectively characterizes solutions of the
nonhomogeneous system \eqref{3.1} described in Theorem \ref{t2.20}.
\begin{theorem}\label{t3.4}
Assume that the Hamiltonian system \eqref{HSa} which satisfies
Hypothesis \ref{h2.4} is either in the limit circle case or in the limit
point case at $-\infty$ and is either in the limit circle or in the limit
point case at $\infty$.  Let $f\in\ell_A^2(\bbZ)$ and let
$y(z,\cdot)$ represent the $\ell_A^2(\bbZ)$-solution of the nonhomogeneous
system \eqref{3.1} which is defined by \eqref{3.30}. Then, $y(z,\cdot)$
represents the unique $\ell_A^2(\bbZ)$-solution of \eqref{3.1} which
satisfies the boundary conditions given by
\begin{align}
\lim_{k\to \infty} \hatt{U}_+^\circledast(z,k) J_\rho \hatt{y}(z,k)&=0,
\label{3.57} \\
\lim_{k\to -\infty} \hatt{U}_-^\circledast(z,k) J_\rho \hatt{y}(z,k)&=0, 
\label{3.58}
\end{align}
where $U_\pm(z,k)=U_\pm(z,k,k_0,\ta)$ and $U_\pm^\circledast(z,k)$ are 
defined in Theorem \ref{t2.20}. Moreover, when \eqref{HSa} is in the limit
point case at $\infty$ $($resp., $-\infty$$)$, the corresponding
boundary condition \eqref{3.57} $($resp., \eqref{3.58}$)$ is superfluous
and can be dropped.
\end{theorem}
\begin{proof}
There are four cases to be considered. However, by symmetry this may be
reduced to three.  The first of these cases we consider assumes that the
limit circle case holds at $-\infty$ while the limit point case holds at
$\infty$.

Let $u(z,k)$ and $v(z,k)$ be $\ell_A^2(\bbZ)$-solutions of \eqref{3.1}
which satisfy \eqref{3.57} and \eqref{3.58}. Then, $w(z,k)=
u(z,k)-v(z,k)$ satisfies \eqref{3.57} and is an $\ell_A^2(\bbZ)$-solution
of the Hamiltonian system \eqref{HSa}, and as a result,
$w(z,k)=U_+(z,k)\zeta$ for some $\zeta\in\bbC^{2m}$. Thus, by
Lemma \ref{l3.2}, \eqref{3.57}, and
\eqref{2.83}, we see that
\begin{align}
0&=\lim_{k\to -\infty} \hatt{U}_-^\circledast(z,k) J_\rho \hatt{w}(z,k)\\
&=\lim_{k\to -\infty} \hatt{U}_-^\circledast(z,k) J_\rho\hatt{U}_+(z,k)\zeta\\
&=[M_-(z) - M_+(z)]\zeta.
\end{align}
Given the invertibility of $(M_- - M_+)$, we see that $\zeta=0$. Note that
the boundary condition at $\infty$ given in \eqref{3.58} was not used in
this argument. However, note that it is automatically satisfied by
$y(z,k)$ due to Corollary~\ref{c3.3}.

Suppose that the limit point case holds at both $\infty$ and at $-\infty$. Then, with $u$, $v$, and $w$ as
previously defined, we again see that $w(z,k)=U_+(z,k)\zeta$ for some $\zeta\in\bbC^{2m}$, but that $\zeta=0$
by the reasoning in the previous case. However,  now note that \eqref{3.57} is automatically satisfied by
$y(z,k)$ by Corollary~\ref{c3.3} as it can be restated when the limit point case holds at $-\infty$.

Lastly, we suppose that the limit circle case holds at both $\infty$ and at $-\infty$. Once again, let
$u(z,k)$ and $v(z,k)$ be $\ell_A^2(\bbZ)$-solutions of \eqref{3.1} which
satisfy \eqref{3.57} and \eqref{3.58} and let $w(z,k)= u(z,k)-v(z,k)$. 
Now note that the columns of $U_-(z,k)$ and of $U_+(z,k)$ together form a
basis for all solutions of the Hamiltonian system \eqref{HSa}. Then,
$w(z,k)=U_-(z,k)\zeta + U_+(z,k)\eta$ for some $\zeta, \eta\in\bbC^{2m}$.
Then, by Lemma \ref{l3.2}, \eqref{3.57}, \eqref{2.83}, and
\eqref{2.101},
\begin{align}
0&=\lim_{k\to -\infty}\hatt{U}_-^\circledast(z,k) J_\rho \hatt{w}(z,k)\\
&=\lim_{k\to -\infty}\hatt{U}_-^\circledast(z,k) J_\rho
\begin{bmatrix} \hatt{U}_-(z,k) &  \hatt{U}_+(z,k) \end{bmatrix}\begin{bmatrix}\zeta\\ \eta \end{bmatrix}\\
&=\begin{bmatrix} 0 &  (M_-(z) - M_+(z)) \end{bmatrix}\begin{bmatrix}\zeta
\\
\eta \end{bmatrix}.
\end{align}
Given the invertibility of $(M_- - M_+)$, we see that $\eta=0$.  By
similar reasoning using \eqref{3.58},
\eqref{2.83}, and \eqref{2.101}, we see that $\zeta = 0$.
\end{proof}

In analogy to the treatment in \cite{HS83}, \cite{HS84} in the continuous
context, we will call the kernel $K(z,\cdot,\cdot)$ defined in
\eqref{2.80} the $2m\times 2m$ Green's matrix of the
Hamiltonian system \eqref{HSa} on $\bbZ$ associated with the boundary
conditions \eqref{3.57} (resp., \eqref{3.58}) if \eqref{HSa} is in the
limit circle case at $\infty$ (resp., $-\infty$). If \eqref{HSa} is in
the limit point case at $\infty$ and $-\infty$, $K(z,\cdot,\cdot)$
represents the unique Green's matrix corresponding to \eqref{HSa} on $\bbZ$.

\medskip

Next, we turn to the analogous considerations for half-lines and
start with the righ half-line $[k_0,\infty)$. We assume Hypothesis
\ref{h2.4} on $[k_0+1,\infty)$ and again consider the nonhomogeneous
system \eqref{3.1} associated with the Hamiltonian system \eqref{HSa}
which is in the limit point or the limit circle case at $\infty$. We
assume that $f(k)$ is defined for $k\in [k_0+1,\infty)$ and that
$f\in\ell_A^2([k_0+1,\infty))$.

We describe for $k,\ell\in [k_0,\infty)$, and $z\in\bbC\backslash\bbR$, a
matrix
$K_+(z,k,\ell)\in\bbC^{2m\times 2m}$ for which the following properties
hold:
\begin{equation}\label{3.65}
\sum_{\ell\in [k_0,\infty)}K_+(z,k,\ell)A(\ell)K_+(z,k,\ell)^*<\infty, 
\quad k\in [k_0,\infty).
\end{equation}
If $f\in\ell_A^2([k_0+1,\infty))$ and if
\begin{equation}\label{3.66}
y(z,k)=\sum_{\ell\in [k_0+1,\infty)} K_+(z,k,\ell)A(\ell)f(\ell), \quad 
k\in [k_0,\infty),
\end{equation}
then $y(z,\cdot)\in\ell_A^2([k_0,\infty))$ and $y(z,\cdot)$ satisfies
\eqref{3.1} on $[k_0+1,\infty)$. In addition, it will be seen that
$y(z,\cdot)$ satisfies certain boundary conditions at $k=k_0$ and at
$\infty$ (if \eqref{HSa} is in the l.c. case at $\infty$).

As in Theorem \ref{t2.20}, we assume for convenience that $z\in\bbC_+$.
\begin{theorem}\label{t2.21}
Assume Hypothesis \ref{h2.4} on $[k_0,\infty)$ and suppose that
$z\in\bbC_+$ and $k,\ell \in[k_0,\infty)$. Let $K_+(z,k,\ell)$ be
defined by
\begin{equation}
K_+(z,k,\ell)=
\begin{cases}
U_+(z,k)\Phi(\bar z,\ell)^*, & k>\ell,\\ 
\begin{pmatrix}
u_{+,1}(z,k)\phi_{1}(\bar z,k)^*& u_{+,1}(z,k)\phi_{2}(\bar z,k)^* \\
\phi_{2}(z,k)u_{+,1}(\bar z,k)^*& \phi_{2}(z,k)u_{+,2}(\bar z,k)^*
\end{pmatrix}, & k=\ell,  \\[1mm]
\Phi(z,k)U_+(\bar z,\ell)^*, & k<\ell. \lb{3.67}
\end{cases}
\end{equation}
Here $U_{+}(z,k)$ is defined in \eqref{2.19} with
$M=M_{+}(z)\in\partial\calD_{+}(z,k_0,\ta)$,
$U_{+}(\bar z,k)$ is defined in \eqref{2.19} with
$M=M_+(\bar z)=M_{+}(z)^*\in\partial\calD_{+}(\bar z,k_0,\ta)$, and
\begin{equation}
U_{+}(z,k)= \begin{pmatrix}u_{+,1}(z,k)\\ u_{+,2}(z,k) \end{pmatrix}
\end{equation}
with $u_{+,j}(z,k)\in\bbC^{m\times m}$, $j=1,2$, $k\in [k_0,\infty)$. In
addition, 
$\Phi(z,k)$ is defined in \eqref{FS}.
With $K_+(z,k,\ell)$ so defined, \eqref{3.65} is satisfied.
Moreover, as defined in \eqref{3.66}, $y(z,\cdot)$
satisfies \eqref{3.1} on $[k_0+1,\infty)$ and is in 
$\ell_A^2([k_0,\infty))$.
\end{theorem}
\begin{proof}
This result follows using the same steps already given for the proof of
Theorem \ref{t2.20} with the following identifications replacing those of
\eqref{2.86}:
\begin{alignat}{2}
\varphi_+(z,k)&=u_{+,1}(z,k), &\qquad \varphi_-(z,k)=\phi_1(z,k),\\
\vartheta_+(z,k)&=u_{+,2}(z,k),&\qquad \vartheta_-(z,k)=\phi_2(z,k).
\end{alignat}
This  assigns the same meaning to $U_{+}(z,k)$ as in Theorem \ref{t2.20},  but
unlike Theorem \ref{t2.20} it makes the further assignment given by
$U_{-}(z,k)=\Phi(z,k)$.

The principal effect of this set of assignments in modifying the proof of
Theorem \ref{t2.20} comes with the
realization that \eqref{2.83} is now replaced by
\begin{align}
\hatt{U}_{-}(\bar z, k)^*J_{\rho}(k)\hatt{U}_{+}(z,k)&=
\hatt{U}_{-}(\bar z, k_0)^*J_{\rho}(k_0)\hatt{U}_{+}(z,k_0)\\
&=\hatt{\Phi}(\bar z, k_0)^*J_{\rho}(k_0)\hatt{U}_{+}(z,k_0)\\
&=I_m,
\intertext{together with}
\hatt{U}_{+}(\bar z, k)^*J_{\rho}(k)\hatt{U}_{-}(z,k)&=
\hatt{U}_{+}(\bar z, k_0)^*J_{\rho}(k_0)\hatt{U}_{-}(z,k_0)\\
&=\hatt{U}_{+}(\bar z, k_0)^*J_{\rho}(k_0)\hatt{\Phi}(z,k_0)\\
&=-I_m.
\end{align}
As a consequence of the identifications now given, we make the
further assignment and modification to the proof
given in the previous theorem: $\omega = I_m$.

As in Theorem\ref{t2.20}, we note that \eqref{2.88} and \eqref{2.89}
imply that $K_+(z,k,\ell)$ satisfies \eqref{3.65}. And finally, by an 
argument in direct analogy
with that given in \cite[Lemma 4.2]{HS83}, or in \cite[Lemma 2.1]{KR74} 
for the one singular
endpoint case,
we see that $y$ given by \eqref{3.66} satisfies
\begin{equation}
\sum_{k\in[k_0+1,\infty)}y(z,k)^*A(k)y(z,k) \le 
(\Im (z))^{-2}\sum_{k\in[k_0+1,\infty)}f(k)^*A(k)f(k).
\end{equation}
As a result, we note that $y(z,\cdot)\in\ell^2_A([k_0,\infty))$ whenever 
$f\in\ell^2_A([k_0+1,\infty))$.
\end{proof}

We now  state a result whose proof is analogous to that of  
Theorem \ref{t3.4}.

\begin{theorem}\label{t3.6}
Assume that the Hamiltonian system \eqref{HSa} which satisfies 
Hypothesis \ref{h2.4} on $[k_0,\infty)$
is in the limit point or limit circle case at $\infty$. Let
$f(k)\in\bbC^{2m}$ be defined for $k\in [k_0+1,\infty)$ with
$f\in\ell_A^2([k_0+1,\infty))$ and let $y(z,k)$
be described by \eqref{3.66}. Then, $y(z,\cdot)\in\ell_A^2([k_0,\infty))$
and $y(z,\cdot)$ satisfies \eqref{3.1} on $[k_0+1,\infty)$.
Moreover, $y(z,\cdot)$ is uniquely defined  by the boundary conditions
\begin{align}
\ta \hatt{y}(z,k_0) &=0, \lb{3.77} \\
\lim_{k\to \infty} \hatt{U}_+^\circledast(z,k) J_\rho \hatt{y}(z,k)&=0, 
\lb{3.78}
\end{align}
where $U_+(z,k)=U_+(z,k,k_0,\ta)$ and $U_+^\circledast(z,k)$ are defined 
in Theorem \ref{t2.20}. When \eqref{HSa} is in the 
limit point case at $\infty$, the corresponding boundary condition in
\eqref{3.78} is superfluous and can be dropped.
\end{theorem}

Again, in analogy to the treatment in \cite{HS81} in the continuous
context, we will call the kernel $K_+(z,\cdot,\cdot)$ defined in
\eqref{3.67} the $2m\times 2m$ half-line Green's matrix of the Hamiltonian
system \eqref{HSa} on $[k_0,\infty)$ associated with the boundary
conditions \eqref{3.77} and \eqref{3.78} (if \eqref{HSa} is in the limit
circle case at $\infty$). 

\medskip

Finally, we briefly turn to the left half-line case $(-\infty,k_0]$. We 
assume Hypothesis \ref{h2.4} on $(-\infty,k_0]$ and again consider the
nonhomogeneous system \eqref{3.1} associated with the Hamiltonian system
\eqref{HSa} which is in the limit point or the limit circle case at
$-\infty$. We assume that $f(k)$ is defined for $k\in (-\infty,k_0-1]$
and that $f\in\ell_A^2((-\infty,k_0-1])$.

We describe for $k,\ell\in (-\infty,k_0]$, and $z\in\bbC\backslash\bbR$,
a matrix
$K_-(z,k,\ell)\in\bbC^{2m\times 2m}$ for which the following properties
hold:
\begin{equation}\label{3.79}
\sum_{\ell\in (-\infty,k_0]}K_-(z,k,\ell)A(\ell)K_-(z,k,\ell)^*<\infty, 
\quad k\in (-\infty,k_0].
\end{equation}
If $f\in\ell_A^2((-\infty,k_0-1])$ and if
\begin{equation}\label{3.80}
y(z,k)=\sum_{\ell\in (-\infty,k_0-1]} K_-(z,k,\ell)A(\ell)f(\ell), \quad 
k\in (-\infty,k_0], 
\end{equation}
then $y(z,\cdot)\in\ell_A^2((-\infty,k_0])$  and $y(z,\cdot)$ satisfies
\eqref{3.1} on $(-\infty,k_0-1]$. In addition, it will be seen that
$y(z,\cdot)$ satisfies certain boundary conditions at $k=k_0$ and at
$-\infty$ (if \eqref{HSa} is in the l.c. case at $-\infty$).

As in Theorems \ref{t2.20} and \ref{t2.21}, we assume for convenience that
$z\in\bbC_+$.
\begin{theorem}
Assume Hypothesis \ref{h2.4} on $(-\infty,k_0]$ and suppose that
$z\in\bbC_+$ and $k,\ell\in (-\infty,k_0]$. Let $K_-(z,k,\ell)$ be
defined by
\begin{equation}
K_-(z,k,\ell)=
\begin{cases}
-\Phi(z,k)U_-(\bar z,\ell)^*, & k>\ell, \\
\begin{pmatrix}
-\phi_{1}(z,k)u_{-,1}(\bar z,k)^*& -\phi_{1}(z,k)u_{-,2}(\bar z,k)^* \\
-u_{-,2}(z,k)\phi_{1}(\bar z,k)^*& -u_{-,2}(z,k)\phi_{2}(\bar z,k)^*
\end{pmatrix}, & k=\ell,  \\[1mm] 
-U_-(z,k)\Phi(\bar z,\ell)^*, & k<\ell. \lb{3.81}
\end{cases}
\end{equation}
Here $U_{-}(z,k)$ is defined in \eqref{2.19} with
$M=M_{-}(z)\in\partial\calD_{-}(z,k_0,\ta)$,
$U_{-}(\bar z,k)$ is defined in \eqref{2.19} with
$M=M_-(\bar z)=M_{-}(z)^*\in\partial\calD_{-}(\bar z,k_0,\ta)$, and
\begin{equation}
U_{-}(z,k)= \begin{pmatrix}u_{-,1}(z,k)\\ u_{-,2}(z,k) \end{pmatrix}
\end{equation}
with $u_{-,j}(z,k)\in\bbC^{m\times m}$, $j=1,2$, $k\in (-\infty,k_0]$, and
$z\in\bbC_+$. In addition, $\Phi(z,k)$ is defined in \eqref{FS}.
With $K_-(z,k,\ell)$ so defined, \eqref{3.79} is satisfied.
Moreover,  as defined in \eqref{3.80}, $y(z,\cdot)$ satisfies \eqref{3.1}
on $(-\infty,k_0-1]$ and is in $\ell_A^2((-\infty,k_0])$.
\end{theorem}
\begin{proof}
As seen in the proof of Theorem \ref{t2.21},
this result follows using the same steps already given for the proof of
Theorem \ref{t2.20} with the following identifications replacing those of
\eqref{2.86}:
\begin{alignat}{2}
\varphi_+(z,k)&=\phi_1(z,k), &\qquad \varphi_-(z,k)=u_{-,1}(z,k),\\
\vartheta_+(z,k)&=\phi_2(z,k), &\qquad \vartheta_-(z,k)=u_{-,2}(z,k).
\end{alignat}
This  assigns the same meaning to $U_{-}(z,k)$ as in Theorem \ref{t2.20},
but unlike Theorem \ref{t2.20} it makes the further assignment given by
$U_{+}(z,k)=\Phi(z,k)$.  As in Theorem \ref{t2.21}, we again find that
\begin{align}
\hatt{U}_{\mp}(\bar z, k)^*J_{\rho}(k)\hatt{U}_{\pm}(z,k)&=
\hatt{U}_{\mp}(\bar z, k_0)^*J_{\rho}(k_0)\hatt{U}_{\pm}(z,k_0)\\
&=\mp I_m,
\end{align}
and hence that  $\omega=-I_m$.

As in Theorem \ref{t2.21}, we note that \eqref{2.88} and \eqref{2.89}
imply that $K_-(z,k,\ell)$ satisfies \eqref{3.79}. And finally, by an
argument in direct analogy with that given  in \cite[Lemma 2.1]{KR74} for
the one singular endpoint case, we see that $y$ given by \eqref{3.80}
satisfies
\begin{equation}
\sum_{k\in (-\infty,k_0-1]}y(z,k)^*A(k)y(z,k) \le (\Im
(z))^{-2}\sum_{k\in (-\infty,k_0-1]}f(k)^*A(k)f(k).
\end{equation}
As a result, we note that $y(z,\cdot)\in\ell^2_A((-\infty,k_0])$ whenever 
$f\in\ell^2_A((-\infty,k_0-1])$.
\end{proof}

Lastly, we  state a result whose proof is again analogous to that of
Theorem \ref{t3.4}.

\begin{theorem}
Assume that the Hamiltonian system \eqref{HSa} which satisfies 
Hypothesis \ref{h2.4} on $(-\infty,k_0]$ is in the limit point or limit
circle case at $-\infty$. Let $f(k)\in\bbC^{2m}$ be defined  for
$k\in (-\infty,k_0-1]$ with $f\in\ell_A^2((-\infty,k_0-1])$ and let $y(z,k)$
be described by \eqref{3.80}. Then, $y(z,\cdot)\in\ell_A^2((-\infty,k_0])$
and $y(z,\cdot)$ satisfies \eqref{3.1} on $(-\infty,k_0-1]$. Moreover,
$y(z,\cdot)$ is uniquely defined  by the boundary conditions
\begin{align}
\ta \hatt{y}(z,k_0) &=0, \lb{3.87} \\
\lim_{k\to -\infty} \hatt{U}_-^\circledast(z,k) J_\rho \hatt{y}(z,k)&=0, 
\lb{3.88}
\end{align}
where $U_-(z,k)=U_-(z,k,k_0,\ta)$ and $U_-^\circledast(z,k)$ are defined 
in Theorem \ref{t2.20}. When \eqref{HSa} is in the limit point case at
$-\infty$, the corresponding boundary condition \eqref{3.88} is
superfluous and can be dropped.
\end{theorem}

As in the previous half-line case, we will call the kernel
$K_-(z,\cdot,\cdot)$ defined in \eqref{3.81} the
$2m\times 2m$ half-line Green's matrix of the Hamiltonian system
\eqref{HSa} on $(-\infty,k_0]$ associated with the boundary conditions
\eqref{3.87} and \eqref{3.88} (if \eqref{HSa} is in the limit circle case
at $-\infty$). 

In our subsequent paper \cite{CGR04}, the explicit formulas \eqref{2.80} 
for the Green's function on $\bbZ$ together with their
asymptotic expnsions as $|z|\to\infty$ will be used to prove trace
formulas of the matrix-valued Jacobi operator \eqref{HSd}, \eqref{2.16}
and the Dirac-type difference expression \eqref{HSd}, \eqref{2.15}. This
in turn then yields Borg-type uniqueness theorems for these Jacobi and
Dirac-type difference operators in analogy to our treatment of
Schr\"odinger and Dirac-type differential operators in \cite{CG02} and
\cite{CGHL00}. As indicated at the end of the introduction, these results
are relevant in connection with the nonabelian Toda and Kac--van
Moerbeke  hierarchies of completely integrable evolution equations.

\medskip
\noindent {\bf Acknowledgements.} We thank Don Hinton and Walter Renger
for many useful discussions. Moreover, we are indebted to Don Hinton for
pertinent hints to the literature.



\begin{thebibliography}{99}
%
\bi{AG94} D.~Alpay and I.~Gohberg, {\it Inverse spectral problems
for difference operators with rational scattering matrix function},
Integr. Equat. Oper. Th. {\bf 20}, 125--170 (1994).
%
\bi{AN84} A.~I.~Aptekarev and E.~M.~Nikishin, {\it The scattering
problem for a discrete Sturm-Liouville problem}, Math. USSR
Sborn. {\bf 49}, 325--355 (1984).
%
\bi{AD56} N.~Aronszajn and W.~F.~Donoghue, {\it On
exponential representations
of analytic functions in the upper half-plane with
positive imaginary
part}, J. Analyse Math. {\bf 5}, 321-388 (1956-57).
%
\bi{As66} T.~Asahi, {\it Spectral theory of the difference equations},
Progr. Theort. Phys. Suppl. {\bf 36}, 55--96 (1966).
%
\bi{At64} F.~V.~Atkinson, {\it Discrete and Continuous
Boundary Problems}, Academic Press, New York, 1964.
%
\bi{At88} F.~V.~Atkinson, {\it Asymptotics of the
Titchmarsh--Weyl function in the matrix case}, unpublished manuscript.
%
\bi{At88a} F.~V.~Atkinson, {\it On the order of magnitude of
Titchmarsh--Weyl functions}, Diff. Integral Eqs. {\bf 1}, 79--96 (1988).
%
\bi{BGMS03} E.\ D.\ Belokolos, F.\ Gesztesy, K.\ A.\ Makarov, and L.\ A.\
Sakhnovich, {\it Matrix-valued generalizations of the theorems of Borg 
and Hochstadt}, in {\it Evolution Equations}, G.\ Ruiz Goldstein, R.\
Nagel, and S.\ Romanelli (eds.), Lecture Notes in Pure and Applied
Mathematics, Vol.\ 234, Marcel Dekker, New York, 2003, p.\ 1--34. 
%
\bi{Be68} Ju.~Berezanskii, {\it Expansions in Eigenfunctions
of Selfadjoint Operators}, Transl. Math. Mongraphs, Vol. 17,
Amer. Math. Soc., Providence, R.I., 1968.
%
\bi{BGS86} Yu.~M.~Berezanskii, M.~I.~Gekhtman, and
M.~E.~Shmoish, {\it Integration of some chains of nonlinear
difference equations by the method of the inverse spectral
problem}, Ukrain. Math. J. {\bf 38}, 74--78 (1986).
%
%
\bi{BG90} Yu.~M.~Berezanskii and M.~I.~Gekhtman, {Inverse problem
of the spectral analysis and non-abelian chains of nonlinear
equations}, Ukrain. Math. J. {\bf 42}, 645--658 (1990).
%
\bi{BGHT98}W.~Bulla, F.~Gesztsy, H.~Holden, and G.~Teschl, {\it Algebro-
Geometric Quasi-Periodic Finte-Gap Solutions of the Toda and Kac--van
Moerbeke Hierarachies,} Memoirs of the Amer. Math. Soc. {\bf 135/641},
(1998).
%
\bi{Ca76} R. W. Carey, {\it A unitary invariant for pairs
of self-adjoint
operators}, J. reine angew. Math. {\bf 283}, 294--312 (1976).
%
\bi{CE67} J.~Chaudhuri and W.~N.~Everitt, {\it On the spectrum of
ordinary second order differential operators}, Proc. Roy. Soc. Edinburgh
{\bf 68A}, 95--119 (1967--68).
%
\bi{CG01} S.~Clark and F.~Gesztesy, {\it Weyl--Titchmarsh
$M$-function asymptotics for matrix-valued Schr\"odinger
operators}, Proc. London Math. Soc. {\bf 82}, 701--724 (2001).
%
\bi{CG02} S.~Clark and F.~Gesztesy, {\it Weyl--Titchmarsh
$M$-function asymptotics and Borg-type theorems for Dirac
operators}, Trans. Amer. Math. Soc. {\bf 354}, 3475--3534 (2002).
%
\bi{CG03} S.~Clark and F.~Gesztesy, {\it On
Povzner--Wienholtz-type Self-Adjointness Results for Matrix-Valued
Sturm--Liouville Operators}, Proc. Roy. Soc. Edinburgh {\bf A} (to
appear).
%
\bi{CGHL00} S.~Clark, F.~Gesztesy, H.~Holden, and
B.~M.~Levitan, {\it Borg-type theorems for matrix-valued
Schr\"odinger operators}, J. Diff. Eqs. {\bf 167}, 181--210 (2000).
%
\bi{CGR04} S.~Clark, F.~Gesztesy, and W.~Renger, {\it Trace formulas and
Borg-type theorems for finite difference operators}, in preparation.
%
\bi{DL96} A.~J.~Duran and P.~Lopez-Rodriguez, {\it Orthogonal
matrix polynomials: zeros and Blumenthal's theorem}, J. Approx. Th.
{\bf 84}, 96--118 (1996).
%
\bi{DL00} A.~J.~Duran and P.~Lopez-Rodriguez, {\it $N$-extremal
matrices of measures for an indeterminate matrix moment problem}, J.
Funct. Anal. {\bf 174}, 301--321 (2000).
%
\bi{DV95} A.~J.~Duran and W.~ Van Assche, {\it Orthogonal matrix
polynomials and higher-order recurrence relations}, Lin. Algebra
Appl. {\bf 219}, 261--280 (1995).
%
\bibitem{Ev59} W.~N.~Everitt, {\it Integrable-square solutions
of ordinary differential equations},
Quart. J. Math. Oxford (2) {\bf 10}, 145--155 (1959).
%
\bibitem{Ev62} W.~N.~Everitt, {\it Integrable-square solutions
of ordinary differential equations II},
Quart. J. Math. Oxford (2) {\bf 13}, 217--220 (1962).
%
\bibitem{Ev63} W.~N.~Everitt, {\it Integrable-square solutions
of ordinary differential equations III},
Quart. J. Math. Oxford (2) {\bf 14}, 170--180 (1963).
%
\bibitem{Ev63a} W.~N.~Everitt, {\it Fourth order singular differential
operators}, Math. Ann. {\bf 149}, 320--340 (1963).
%
\bibitem{Ev64} W.~N.~Everitt, {\it Singular differential
equations I: The even order case}, Math. Ann. {\bf 156}, 9-24 (1964).
%
\bibitem{Ev67} W.~N.~Everitt, {\it Singular differential
equations II: Some self-adjoint even order cases}, Quart. J. Math. Oxford
(2) {\bf 18}, 13--32 (1967).
%
\bi{Ev72} W.~N.~Everitt, {\it On a property of the $m$-coefficient
of a second-order linear differential equation}, J. London Math. 
Soc. (2), {\bf 4}, 443--457 (1972).
%
\bibitem{Ev72a} W.~N.~Everitt, {\it Integrable-square, analytic solutions
of odd-order, formally symmetric, ordinary differential equations},
Proc. London Math. Soc. (3) {\bf 25}, 156--182 (1972).
%
\bi{EB80} W.~N.~Everitt and C.~Bennewitz, {\it Some remarks on the
Titchmarsh--Weyl $m$-coefficient}, in {\it Tribute to {\AA}ke Pleijel},
Mathematics Department, University of Uppsala, Sweden, 1980, pp.\ 49--108.
%
\bi{EH78} W.~N.~Everitt and S.~G.~Halvorsen, {\it On the asymptotic
form of the Titchmarsh--Weyl $m$-coefficient}, Appl. Anal. {\bf 8},
153--169 (1978).
%
\bi{EHS83} W.~N.~Everitt, D.~B.~Hinton, and J.~K.~Shaw, {\it The
asymptotic form of the Titchmarsh--Weyl coefficient for Dirac 
systems}, J. London Math. Soc. (2), {\bf 27}, 465--476 (1983).
%
\bibitem{EK76} W.~N.~Everitt and K.~Kumar, {\it On the Titchmarsh--Weyl
theory of ordinary symmetric differential expressions. I. The odd-order
case}, Nieuw Arch. Wisk. (3) {\bf 24}, 109--145 (1976).  
%
\bibitem{EK76a} W.~N.~Everitt and K.~Kumar, {\it On the Titchmarsh--Weyl
theory of ordinary symmetric differential expressions. II. The general
theory}, Nieuw Arch. Wisk. (3) {\bf 24}, 1--48 (1976).  
%
\bibitem{EZ79} W.~N.~Everitt and A.~Zettl, {\it Generalized symmetric
ordinary differential expressions. I. The general theory}, Nieuw Arch.
Wisk. (3) {\bf 27}, 363--397 (1979).  
%
\bi{Fu76} M.~Fukushima, {\it A spectral representation on ordinary
linear difference equation with operator-valued coefficients of the
second order}, J. Math. Phys. {\bf 17}, 1084--1072 (1976).
%
\bi{Ge82} J.~S.~Geronimo, {\it Scattering theory and matrix
orthogonal polynomials on the real line}, Circuits Syst. Signal
Process. {\bf 1}, 471--495 (1982).
%
\bi{GH97} F.~Gesztesy and H.~Holden, {\it On trace formulas for
Schr\"odinger-type operators}, in Multiparticle Quantum
Scattering with Applications to Nuclear, Atomic and Molecular
Physics, D.~G.~Truhlar and B.~Simon (eds.), Springer, New
York, 1997, pp.~121--145.
%
\bi{GHSZ93} F.~Gesztesy, H.~Holden, B.~Simon, and Z.~Zhao, {\it On the
Toda and Kac--van Moerbeke systems}, Trans. Amer. Math. Soc. {\bf 339},
849--868 (1993).
%
\bi{GKM02} F.~Gesztesy, A.~Kiselev, and K.~A.~Makarov,
{\it Uniqueness Results for Matrix-Valued Schr\"odinger,
Jacobi, and Dirac-Type Operators}, Math. Nachr. {\bf 239--240},
103--145 (2002).
%
\bi{GKT96} F.~Gesztesy, M.~Krishna, and G.~Teschl, {\it On
isospectral sets of Jacobi operators,} Commun. Math. Phys.
{\bf 181}, 631--645 (1996).
%
\bi{GMT98} F.~Gesztesy, K.~A.~Makarov, and E.~Tsekanovskii,
{\it An Addendum to Krein's Formula,}
J. Math. Anal. Appl. {\bf 222}, 594--606 (1998).
%
\bi{GS03} F.\ Gesztesy and L.\ A.\ Sakhnovich, {\it A class of
matrix-valued Schr\"odinger operators with prescribed finite-band spectra}, 
in {\it Reproducing Kernel Hilbert Spaces, Positivity, System Theory and
Related Topics}, D.\ Alpay (ed.), Operator Theory: Advances and
Applications, Vol.\ 143, Birkh\"auser, Basel, 2003, p.\ 213--253.
%
\bi{GS97} F.~Gesztesy and B.~Simon, {\it{$m$-functions and inverse
spectral analysis for finite and semi-infinite Jacobi matrices}},
J. d'Anal. Math. {\bf 73} (1997), 267--297.
%
\bi{GS00} F.~Gesztesy and B.~Simon, {\it On local
Borg-Marchenko uniqueness results}, Commun. Math. Phys.
{\bf 211}, 273--287 (2000).
%
\bi{GT96} F.~Gesztesy and G.~Teschl, {\it Commutation methods for Jacobi
operators}, J. Diff. Eqs. {\bf 128}, 252--299 (1996).
%
\bi{GT00} F.~Gesztesy and E.~Tsekanovskii, {\it On
matrix-valued Herglotz functions}, Math. Nachr. {\bf 218}, 61--138
(2000).
%
\bi{Hi69} E.~Hille, {\it Lectures on Ordinary Differential
Equations}, Addison-Wesley, Reading, 1969.
%
\bi{HS93}  D.~B.~Hinton and A. Schneider, {\it On the Titchmarsh--Weyl
coefficients for singular S-Hermitian Systems I}, Math.
Nachr. {\bf 163}, 323--342 (1993).
%
\bi{HS97}  D.~B.~Hinton and A. Schneider, {\it On the Titchmarsh--Weyl
coefficients for singular S-Hermitian Systems II}, Math.
Nachr. {\bf 185}, 323--342 (1997).
%
\bi{HS81} D.~B.~Hinton and J.~K.~Shaw, {\it On Titchmarsh--Weyl
$M(\lambda)$-functions for linear Hamiltonian systems},
J. Diff. Eqs. {\bf 40}, 316--342 (1981).
%
\bi{HS82} D.~B.~Hinton and J.~K.~Shaw, {\it On the spectrum
of a singular Hamiltonian system}, Quaest. Math.
{\bf 5}, 29--81 (1982).
%
\bi{HS83} D.~B.~Hinton and J.~K.~Shaw, {\it Hamiltonian
systems of limit point or limit circle type with both endpoints
singular}, J. Diff. Eqs. {\bf 50}, 444--464 (1983).
%
\bi{HS83a} D.~B.~Hinton and J.~K.~Shaw, {\it Parameterization of the
$M(\lambda)$ function for a Hamiltonian system of limit circle type}, 
Proc. Roy. Soc. Edinburgh {\bf 93A}, 349--360 (1983).
%
\bi{HS84} D.~B.~Hinton and J.~K.~Shaw, {\it On boundary
value problems for Hamiltonian systems with two singular points},
SIAM J. Math. Anal. {\bf 15}, 272--286 (1984).
%
\bi{HS86} D.~B.~Hinton and J.~K.~Shaw, {\it On the
spectrum of a singular Hamiltonian system, II}, Quaest. Math. {\bf 10},
1--48 (1986).
%
\bi{JNO00} R.~Johnson, S.~Novo, and R.~Obaya, {\it Ergodic
properties and Weyl $M$-functions for random linear
Hamiltonian systems}, Proc. Roy. Soc. Edinburgh {\bf 130A},
1045--1079 (2000).
%
\bi{KR74} V.~I.~Kogan and F.~S.~Rofe-Beketov, {\it On square-integrable
solutions of symmetric systems of differential equations
of arbitrary order}, Proc. Roy. Soc. Edinburgh {\bf 74A}, 1--40 (1974).
%
\bi{KM98} A.~G.~Kostyuchenko and K.~A.~Mirzoev, {\it Three-term
recurrence relations with matrix coefficients. The completely indefinite
case}, Math. Notes {\bf 63}, 624--630 (1998).
%
\bi{KM99} A.~G.~Kostyuchenko and K.~A.~Mirzoev, {\it Generalized Jacobi
matrices and deficiency numbers of ordinary differential
operators with polynomial coefficients}, Funct. Anal. Appl. {\bf 33},
25--37 (1999).
%
\bi{KM01} A.~G.~Kostyuchenko and K.~A.~Mirzoev, {\it Complete
indefiniteness tests for Jacobi matrices with matrix entries}, Funct.
Anal. Appl. {\bf 35}, 265--269 (2001).
%
\bi{KS88} S.~Kotani and B.~Simon, {\it Stochastic
Schr\"odinger
operators and Jacobi matrices on the strip}, Commun. Math.
Phys. {\bf 119}, 403--429 (1988).
%
\bi{Kr89a} A.~M.~Krall, {\it $M(\lambda)$ theory for singular Hamiltonian
systems with one singular point}, SIAM J. Math. Anal. {\bf 20},
664--700 (1989).
%
\bi{Kr89b} A.~M.~Krall, {\it $M(\lambda)$ theory for singular Hamiltonian
systems with two singular points}, SIAM J. Math. Anal. {\bf 20}, 701--715
(1989).
%
\bi{LM00} M.~Lesch and M.~Malamud, {\it The inverse spectral
problem for first order systems on the half line}, Operator
Theory: Advances and Applications, Vol.~117, Birkh\"auser,
Basel, 2000, pp.~199--238.
%
\bi{LM02} M.~Lesch and M.~Malamud, {\it On the deficiecy indices and
self-adjointness of symmetric Hamiltonian systems}, J. Diff. Eq. {\bf
189}, 556--615 (2003).
%
\bi{Lo99} P.~L{\'o}pez-Rodriguez, {\it Riesz's theorem for
orthogonal matrix polynomials}, Constr. Approx. {\bf 15},
135--151 (1999).
%
\bi{MBO92} F.~G.~Maksudov, E.~M.~Bairamov, and R.~U.~Orudzheva, {\it
The inverse scattering problem for an infinite Jacobi matrix with
operator elements}, Russ. Acad. Sci. Dokl. Math. {\bf 45},
366--370 (1992).
%
\bi{Or76} S.~A.~Orlov, {\it Nested matrix disks analytically depending
on a parameter, and theorems on the invariance of ranks of radii of
limiting disks}, Math. USSR Izv. {\bf 10}, 565--613 (1976).
%
\bi{Os97} A.~S.~Osipov, {\it Integration of non-abelian Lanmuir type
lattices by the inverse spectral problem method}, Funct. Anal. Appl. {\bf
31}, 67--70 (1997). 
%
\bi{Os00} A.~S.~Osipov, {\it Some properties of resolvent sets of
second-order difference operators with matrix coefficients}, Math. Notes
{\bf 68}, 806--809 (2000). 
%
\bi{Os02} A.~Osipov, {\it On some issues related to the moment problem
for the band matrices with operator elements}, J. Math. Anal. Appl. {\bf
275}, 657--675 (2002).
%
\bi{Ro60} F.~S.~Rofe-Beketov, {\it Expansions in eigenfunctions of
infinite systems of differential equations in the non-self-adjoint
and self-adjoint cases}, Mat. Sb. {\bf 51}, 293--342 (1960).
(Russian.)
%
\bi{Sa92} A.~L.~Sakhnovich, {\it Spectral functions of a
canonical  system of order $2n$}, Math. USSR Sbornik {\bf 71},
355--369  (1992).
%
\bi{Sa02} A.~Sakhnovich, {\it Dirac type and canonical systems: spectral,
and Weyl--Titchmarsh functions, direct and inverse problems}, Inverse
Problems {\bf 18}, 331--348 (2002).
%
\bi{Sa94a}  L.~A.~Sakhnovich, {\it Method of operator identities
and problems of analysis}, St. Petersburg Math. J. {\bf 5}, 1--69
(1994).
%
\bi{Sa97}  L.~A.~Sakhnovich, {\it Interpolation Theory and its
Applications}, Kluwer, Dordrecht, 1997.
%
\bi{Sa99a}  L.~A.~Sakhnovich, {\it Spectral Theory of Canonical
Differential Systems. Method of Operator Identities}, Operator
Theory: Advances and Applications, Vol. 107, Birkh\"auser, Basel,
1999.
%
\bi{Se80} V.~P.~Serebrjakov, {\it The inverse problem of scattering
theory for difference equations with matrix coefficients}, Sov. Math.
Dokl. {\bf 21}, 148--151 (1980).
%
\bi{Te00} G.~Teschl, {\it Jacobi Operators and Completely Integrable
Nonlinear Lattices}, Mathematical Surveys and Monographs, {\b 72},
American Mathematical Society, (2000).
%
\bi{To89} M. Toda, {\it Theory of Nonlinerar Lattices}, 2nd
enl. ed., Springer, Berlin, 1989.
%
\bi{We87} J.~Weidmann, {\it Spectral Theory of Ordinary Differential
Operators}, Lecture Notes in Math. {\bf 1258}, Springer, Berlin,
1987.
%
\bi{We10} W.~Weyl, {\it \"Uber gew\"ohnliche Differentialgleichungen mit
Singularit\"aten und die zugeh\"origen Entwicklungen willk\"urlicher
Funktionen}, Math. Ann. {\bf 68}, 220--269  (1910).
%
\end{thebibliography}
\end{document}